\documentclass[11pt]{article}
\usepackage{amsmath,amssymb,euscript}
\usepackage{graphicx}
\usepackage[applemac]{inputenc}

\newcommand{\ol}[1]{\overline{#1}}

\newcommand{\ds}{\displaystyle}
\newcommand{\rk}[1]{\mathsf{rank}\left ( #1 \right )}
\renewcommand{\span}[1]{\mathsf{span}\left \{ #1 \right \}}
\newcommand{\diag}[1] {\mathsf{diag}\{ #1 \}}
\newcommand{\rsm}[1]{\mathsf{R-Smith}\left( #1\right)}
\newcommand{\lsm}[1]{\mathsf{L-Smith}\left( #1\right)}
\newcommand{\ddt}{\frac{d}{dt}}

\newcommand{\dg}[1]{{\mathfrak d}\left(#1\right)}
\newcommand{\dgot}{{\mathfrak d}}
\newcommand{\lin}[2]{{\cal L}_{#1}\left( #2 \right)}
\newcommand{\ord}[1]{\mathsf{ord}( #1 )}
\newcommand{\CC}{{\mathbb C}}

\newcommand{\II}{{\EuScript I}}

\newcommand{\kk}{{\mathfrak K}}
\newcommand{\MM}{{\EuScript M}}

\newcommand{\NN}{{\mathbb N}}

\newcommand{\RR}{{\mathbb R}}

\newcommand{\T}{{\mathrm T}}

\newcommand{\TTT}{{\mathfrak T}}
\newcommand{\UU}{{\EuScript U}}

\newcommand{\XX}{{\mathfrak X}}
\newcommand{\XXX}{{\EuScript X}}
\newcommand{\YY}{{\mathfrak Y}}
\newcommand{\YYY}{{\EuScript Y}}

\newcommand{\ZZZ}{{\EuScript Z}}

\newtheorem{thm}{Theorem}
\newtheorem{pr}{Proposition}
\newtheorem{cor}{Corollary}
\newtheorem{lem}{Lemma}
\newtheorem{defn}{Definition}
\newtheorem{rem}{Remark}

\title{On Necessary and Sufficient Conditions for Differential Flatness}
\author{
Jean L{\'e}vine 
            \thanks{
            Centre Automatique et Syst\`{e}mes, Unit\'{e} Math\'{e}matiques et Syst\`{e}mes,
            MINES-ParisTech, 35  rue Saint--Honor{\'e}, 77300 Fontainebleau, France.
              Tel.: +33 1 64 69 48 58.
              E-mail: {\tt{jean.levine@mines-paristech.fr}}
              }
}

\date{To appear in Applicable Algebra in Engineering, Communication and Computing,\\January 2011.} 

\begin{document}
\maketitle

\begin{abstract}
This paper is devoted to the characterization of differentially flat nonlinear systems in implicit representation, after elimination of the input variables, in the differential geometric framework of manifolds of jets of infinite order.  
We extend the notion of Lie-B\"acklund equivalence, introduced in \cite{FLMR-ieee}, to this implicit context and focus attention on Lie-B\"acklund isomorphisms associated to flat systems, called \emph{trivializations}. They can be locally characterized in terms of polynomial matrices of the indeterminate $\ddt$, whose range is equal to the kernel of the polynomial matrix associated to the implicit variational system. Such polynomial matrices are useful to compute the ideal of differential forms generated by the differentials of all possible trivializations. We introduce the notion of a \emph{strongly closed ideal} of differential forms, and prove that flatness is equivalent to the strong closedness of the latter ideal, which, in turn, is equivalent to the existence of solutions of the so-called \emph{generalized moving frame structure equations}. Two sequential procedures to effectively compute flat outputs are deduced and various examples and consequences are presented. 
\end{abstract}


\section*{Introduction}
Differential flatness, or more shortly, flatness, is a system property introduced more than ten years ago in \cite{pM1,FLMR1}. 

Let us briefly state an informal definition, that will be made more precise later: given a nonlinear system 
\begin{equation}\label{sys0}
\dot{x}=f(x,u)
\end{equation}
where $x=(x_{1},\ldots, x_{n})$ is the state (belonging to a given smooth $n$-dimensional manifold) and $u=(u_{1},\ldots, u_{m})$ the control vector, $m\leq n$, the system (\ref{sys0}) is said to be locally (differentially) flat if, and only if, there exists a vector
$y=(y_{1},\ldots,y_{m})$ such that
\begin{itemize}
\item $y$ and its successive time derivatives $\dot{y},\ddot{y},\ldots,$ are locally independent,
\item $y$ is locally a function of $x$, $u$ and a finite number of time derivatives of the components of $u$,
\item $x$ and $u$ can be locally expressed as functions of the components of $y$ and a finite number of their derivatives: $x=\varphi_{0}(y,\dot{y},\ldots, y^{(\alpha)})$, $u=\varphi_{1}(y,\dot{y},\ldots, y^{(\alpha+1)})$, for some multi-integer $\alpha=(\alpha_{1},\ldots,\alpha_{m})$, and with the notation $y^{(\alpha)}=(\frac{d^{\alpha_{1}} y_{1}}{dt^{\alpha_{1}}},\ldots, \frac{d^{\alpha_{m}}y_{m}}{dt^{\alpha_{m}}})$.
\end{itemize}
The vector $y$ is called a \emph{flat output}.

This concept has inspired an important literature. See \cite{fliess-levine-martin-ollivier-rouchon-mtns96,MMR-ecc,rudolph-book-1,rudolph-book-2,SiRaAgra,L_book} for surveys on flatness and its applications. To mention just one fact, flatness  provides significant simplifications to the motion planning problem and to several aspects of feedback design. 

Various formalisms have been introduced: finite dimensional differential geometric approaches \cite{CLM,JF,Sh,Sluis-scl93,SchSch}, differential algebra and related approaches \cite{FLMR-ijc,ABMP-ieee,Jak-ecc93}, infinite dimensional differential geometry of jets and prolongations \cite{FLMR-ieee,vNRM-siam,Po-banach93,PdSCF,RM}. 

Note that the notion of flatness may be tied up to two different trends. \\
The first one refers to  the equivalence between underdetermined differential systems, whose archetype is the problem of Monge, their reduction to normal forms and their integration (see e.g. the major contributions \cite{Monge,Serret,Hadam,Goursat,Zerv-1913,H,C,Zerv-1932}). 
The link between this integrability problem and flatness is particularly clear via \'{E}.~Cartan's notion of \emph{absolute equivalence} \cite{C}  (already noted by W. Shadwick \cite{Sh}) and via D. Hilbert's concept of \emph{invertible and without integral} transformations (\emph{umkehrbar, integrallos transformationen} in German) \cite{H}. \\
The second one is related to the notion of parameterization: using the definition presented in this paper, with, in place of (\ref{sys0}), the set of $n-m$ implicit equations (\ref{implsys}) introduced in the next section, where the control variables $u$ are eliminated, flatness may be seen as a generalization in the framework of manifolds of jets of infinite order of the \emph{uniformization of analytic functions} of Hilbert 's 22nd problem \cite{Hilpro}, solved by Poincaré \cite{Poi} in 1907 (see \cite{Bers} for a modern presentation of this subject and recent extensions and results). This problem consists, roughly speaking, given a set of complex polynomial equations in one complex variable, in finding an open dense subset $D$ of the complex plane $\CC$ and a holomorphic function $s$ from $D$ to $\CC$ such that $s$ is surjective and $s(p)$ identically satisfies the given equations for all values of the ``parameter'' $p\in D$.  In our setting, $\CC$ is replaced by a (real) manifold of jets of infinite order, a flat output $y_{1},\ldots,y_{m}$ plays the role of the parameter $p$ and $s$ is the associated Lie-B\"acklund isomorphism $s=(\varphi_{0},\varphi_{1},\dot{\varphi}_{1},\ddot{\varphi}_{1},\ldots)$ with $\varphi_{0}$ and $\varphi_{1}$ defined above.

In the framework of linear finite or infinite dimensional systems, the notions of flatness and  \emph{parameterization} coincide as remarked by \cite{Pom,Pom-Quad}, and in the behavioral approach of \cite{PW}, flat outputs correspond to \emph{latent variables} of \emph{observable image representations} \cite{Trent} (see also \cite{Fl-scl2} for a module theoretic interpretation of the behavioral approach).

The characterization of differentially flat systems has aroused many contributions \cite{ABMP-ieee,CLMscl,CLM,Cht,JF,Jak-ecc93,martin-rouchon-mtns93,PdS-mtns,Po-cocv97,RM,R1,Sh,Sluis-scl93}. Though general necessary and sufficient conditions exist (see e.g. \cite{ABMP-ieee,Cht,PdS-mtns}), they don't provide a practically computable set of conditions. More precisely, \cite{ABMP-ieee} gives an algorithm to compute a basis of the cotangent module, called infinitesimal Brunovsky form, and further integrability conditions are needed to deduce flat outputs. This result has recently been improved by \cite{Cht,Cht-diffiety} using tools from symmetry groups, and by \cite{AP}. In \cite{PdS-mtns} the author proposes to express all the differential relations between the system variables $x,u$ and a candidate flat output $y$ and use Cartan-K\"ahler theory.

We adopt here the formalism of manifolds of jets of infinite order  \cite{AIb,FLMR-ieee,KLV,L_book,Po-banach93,Z} and, as previously mentioned, in place of systems in explicit form (\ref{sys0}), we consider (locally equivalent) implicit systems obtained from (\ref{sys0}) by eliminating the input vector $u$. The main advantage of this representation is to deal with a system described by a smaller number of variables and relations, which significantly reduces the computational burden. We adapt the notions of Lie-B\"acklund equivalence and Lie-B\"acklund isomorphism in this context and show, after restricting to the category of meromorphic functions, that the flatness property is naturally described in terms of polynomial matrices and differential forms deduced from the variational system equations. For a detailed presentation of polynomial rings and non commutative algebra, the reader may refer to \cite{Co,KS} and for exterior differential systems to \cite{BCG3}.

Though our results show some parallelism with those of \cite{ABMP-ieee,Cht}, in particular concerning the study of variational properties, they differ from the latter by the fact that, as previously announced, our computations involve a smaller number of variables, and exploit different ideas, generalizing the linear approach presented in \cite{LNg} (see also \cite{L_book}), and making an extensive use of polynomials of the operator $\ddt$, which turns out to provide more effective flatness conditions as attested by the examples of the last section. 

The paper is organized as follows: the first Section is devoted to the basic description of implicit control systems on manifolds of jets of infinite order. In Section 2, we extend the notions of Lie-B\"acklund equivalence and Lie-B\"acklund isomorphism to the implicit system framework. Section 3 deals with the presentation of some variational properties of flat systems. The necessary and sufficient conditions for flatness are stated in Theorems~\ref{CNSthm} and \ref{movingframe-thm} of Section~4 and some consequences are presented. Section~5 is then devoted to examples and some concluding remarks are given. An appendix on polynomial matrices and their Smith decomposition is provided.


\section{Implicit control systems on manifolds of jets of infinite order}
Given an infinitely differentiable manifold $X$ of dimension $n$, we denote its tangent space at an arbitrary point $x\in X$ by $\T_{x}X$, and its tangent bundle by $\T X = \bigcup_{x\in X} \T_{x}X$ (identified with the vector bundle $\T X \stackrel{{\EuScript P}}{\rightarrow} X$). Let $F$ belong to $C^{\infty}(\T X ;\RR^{n-m})$, the set of  $C^{\infty}$ mappings from $\T X$ to $\RR^{n-m}$. 
In the sequel, we call \textit{smooth} any function of class $C^{\infty}$ in all its variables. 

 We consider an underdetermined implicit system of the form
\begin{equation}\label{implsys}
F(x,\dot{x})=0
\end{equation}
\emph{regular} in the sense that $\rk{\frac{\partial F}{\partial \dot{x}}}=n-m$ in a suitable open subset of $\T X$. 
\begin{rem}\label{implrem}
Note that any explicit system of the form $\dot{x}=f(x,u)$ with $f$ smooth,  $f(x,u)\in \T_{x}X$ for every $x\in X$ and $u$ in an open subset $U$ of $\RR^{m}$, and $\rk{\frac{\partial f}{\partial u}}=m$ in a suitable open subset of $X\times U$, can be locally transformed into (\ref{implsys}):  permuting the lines of $f$, if necessary, such that the $m$ last lines of $f$ are locally independent functions of $u$, and still noting $x$ the permuted vector (this abuse of notations being unambiguous), thanks to the implicit function theorem one gets\mbox{ $u=\mu(x, \dot{x}_{n-m+1},\ldots, \dot{x}_{n})$}, $\mu$ smooth,  and $\dot{x}_{i}=f_{i}(x,\mu(x,\dot{x}_{n-m+1},\ldots, \dot{x}_{n}))$ for $i=1,\ldots, n-m$. Thus, setting $F_{i}(x,\dot{x})=\dot{x}_{i}-f_{i}(x,\mu(x,\dot{x}_{n-m+1},\ldots, \dot{x}_{n}))$, $i=1,\ldots, n-m$, the system is in the implicit form (\ref{implsys}), and since $\frac{\partial F}{\partial \dot{x}}=\diag{I_{n-m},G}$, $G$ being the matrix made of the entries $-\sum_{k=1}^{m} \frac{\partial f_{i}}{\partial u_{k}}\frac{\partial \mu_{k}}{\partial \dot{x}_{j}}$ for $i,j=n-m+1, \ldots, n$, which don't depend on $(\dot{x}_{1},\ldots, \dot{x}_{n-m})$, we have $\rk{\frac{\partial F}{\partial \dot{x}}}=n-m$.

Conversely, (\ref{implsys}) can be transformed into the explicit system $\dot{x}=f(x,u)$ with $f(x,u)\in \T_{x}X$ for every $x\in X$ and  every $u$ in an open subset $U$ of $\RR^{m}$ with $\rk{\frac{\partial f}{\partial u}}=m$:  the condition $\rk{\frac{\partial F}{\partial \dot{x}}}=n-m$ and the implicit function theorem yield $\dot{x}_{i}=f_{i}(x,\dot{x}_{n-m+1},\ldots, \dot{x}_{n})$, \mbox{$i=1,\ldots, n-m$},  with $f_{i}$ smooth for $i=1,\ldots n-m$, and, setting $\dot{x}_{n-m+j}=u_{j}$ for $j=1,\ldots, m$ (thus $u=(u_{1},\ldots,u_{m})\in \RR^{m}$), we finally get $\dot{x}=f(x,u)$ with $f_{n-m+j}(x,u)=u_{j}$ for $j=1,\ldots,m$, and, by definition, $f$ is a smooth vector field on $X$ for every $u$ in an open subset $U$ of $\RR^{m}$, with $\frac{\partial f}{\partial u}=I_{m}$ the identity matrix of $\RR^{m}$. 
\end{rem}
Clearly, the \emph{underdetermined} character of (\ref{implsys}) is expressed by the rank condition $\rk{\frac{\partial F}{\partial \dot{x}}}=n-m$, which means that the system effectively depends on $m$ independent control variables.

We also introduce the following definition:
\begin{defn}
\begin{itemize}
\item A smooth vector field $f$ that depends, for every $x\in X$, on $m$ independent variables $u\in \RR^{m}$ with $\rk{\frac{\partial f}{\partial u}}=m$ in a suitable open set of $X\times \RR^{m}$, is called \emph{compatible} with (\ref{implsys}) if, and only if, it satisfies $F(x,f(x,u))=0$ for every $u\in U$, where $U$ is a suitable open set of $\RR^{m}$.
\item We say that system (\ref{sys0}) admits the local representation (\ref{implsys}) in a neighborhood of $(x_{0},u_{0})$, if, and only if, $f$ is compatible with (\ref{implsys}) in this neighborhood.
\end{itemize}
\end{defn}
In other words, every smooth integral curve of (\ref{sys0}) is a smooth integral curve of (\ref{implsys}) passing through $(x_{0},\dot{x}_{0})$, with $\dot{x}_{0}=f(x_{0},u_{0})$ for every $(x_{0},u_{0})$ in a suitable open set, and vice-versa. \\

In \cite{FLMR-ieee}, infinite systems of coordinates $(x,\ol{u})=(x,u,\dot{u},\ldots)$ have been introduced to deal with prolonged vector fields $\ol{f}(x,\ol{u})= \sum_{i=1}^{n}f_{i}(x,u)\frac{\partial}{\partial x_{i}} + \sum_{j=1}^{m}\sum_{k\geq 0}u_{j}^{(k+1)}\frac{\partial}{\partial u_{j}^{(k)}}$, associated to explicit systems $\dot{x}=f(x,u)$ (see also \cite{Po-banach93} where a similar approach has been developed  independently).

Here, we adopt an external description\footnote{in the sense that the system manifold, namely the set of points $\ol{x}=(x,\dot{x},\ddot{x},\ldots)$ such that $F(x,\dot{x})=0$, is described by means of the larger manifold $\XX =  X\times \RR^{n}_{\infty}$.} of the prolonged manifold containing the solutions of (\ref{implsys}): we consider the infinite dimensional manifold $\XX \triangleq X\times \RR^{n}_{\infty} \triangleq X\times \RR^{n}\times\RR^{n}\times \ldots$ made of the cartesian product of $X$ with a countably infinite number of copies of $\RR^{n}$. To endow $\XX$ with a suitable topology, we define the continuity and differentiability of functions from $\XX$ to $\RR$ as follows:
\begin{defn}\label{contdef}
We say that a function $\varphi$ from $\XX$ to $\RR$ is \emph{continuous} (resp. \emph{differentiable}) if $\varphi$ depends only on a finite (but otherwise arbitrary) number of variables and is continuous (resp. differentiable) with respect to these variables. 
\end{defn}
Thus $\XX$ is endowed with the coarsest topology that makes projections to $X\times \RR^{kn}$ continuous for any $k\in \NN$. This topology can be identified with the infinite product topology of  $X\times \RR^{n}_{\infty}$, each factor being endowed with its natural finite dimensional topology (see e.g. \cite{KLV,Z,L_book}).

$C^{\infty}$ or analytic or meromorphic functions from $\XX$ to $\RR$ are then defined in the usual finite dimensional way since they only depend on a finite (but otherwise arbitrary) number of variables.

We assume that we are given the infinite set of global coordinates of $\XX$:\begin{equation}\label{coord}
\ol{x}\triangleq (x_{1},\ldots,x_{n},\dot{x}_{1},\ldots,\dot{x}_{n},\ddot{x}_{1},\ldots, \ddot{x}_{n}, \ldots, x_{1}^{(k)},\ldots, x_{n}^{(k)},\ldots)
\end{equation} 
and we endow $\XX$ with the so-called trivial Cartan vector field \cite{KLV,Z}
\begin{equation}\label{Cartanvf}
\tau_{\XX}= \sum_{i=1}^{n} \sum_{j\geq 0}x_{i}^{(j+1)}\frac{\partial}{\partial x_{i}^{(j)}}.
\end{equation}
We denote by
$$\frac{d \varphi}{dt} \triangleq  L_{\tau_{\XX}}\varphi=\sum_{i=1}^{n} \sum_{j\geq 0}x_{i}^{(j+1)}\frac{\partial \varphi}{\partial x_{i}^{(j)}}$$ 
the Lie derivative of a function $\varphi\in C^{\infty}(\XX,\RR)$ along the vector field $\tau_{\XX}$ (this series having only a finite number of non zero terms according to Definition~\ref{contdef}).
Therefore, since  $\frac{d}{dt}x_{i}^{(j)} = L_{\tau_{\XX}}x_{i}^{(j)} = \dot{x}_{i}^{(j)}=x_{i}^{(j+1)}$, the Cartan vector field acts on coordinates as a shift to the right. The pair $(\XX, \tau_{\XX})$ is called \emph{manifold of jets of infinite order} or \emph{diffiety} (see \cite{KLV,Z}). For simplicity's sake, we will keep noting $\XX$ alone in place of $(\XX, \tau_{\XX})$ for this manifold.

From now on, $\ol{x}$ (resp. $\ol{y}, \ol{z}, \ldots$)  stands for the sequence of jets of infinite order of $x$ (resp. $y, z, \ldots$).

\begin{defn}\label{sysdef}
A regular implicit control system is defined as a triple $(\XX,\tau_{\XX},F)$ with \mbox{$\XX=X\times \RR^{n}_{\infty}$, $\tau_{\XX}$} the trivial Cartan field on $\XX$, and where $F\in C^{\infty}(\T X ;\RR^{n-m})$ satisfies $\rk{\frac{\partial F}{\partial \dot{x}}}=n-m$ in a suitable open dense subset of $\T X$. 
\end{defn}
Note that if $t\mapsto (x(t),\dot{x}(t))$ is an integral curve of (\ref{implsys}), then $t\mapsto \ol{x}(t)$ is an integral curve of $L_{\tau_{\XX}}^{k}F=0$ for every $k$. Therefore, every integral curve of the system $(\XX,\tau_{\XX},F)$ lies in the set of $\ol{x}$ such that $L_{\tau_{\XX}}^{k}F=0$ for every $k$.

The system $(\RR^{m}_{\infty},\tau_{m},0)$, where $\tau_{m}$ denotes the trivial Cartan field of $\RR^{m}_{\infty}$, is called the \emph{trivial system} of dimension $m$ since it corresponds to the void (unconstrained) implicit system $0=0$. Note that since $n-m=0$, the rank condition on $\frac{\partial F}{\partial \dot{x}}$ which, in this case, is the matrix whose entries are all identically 0, is globally satisfied. The absence of equations relating the coordinates of $\RR^{m}_{\infty}$ makes this system also an explicit one and, in the notations of \cite{FLMR-ieee}, it corresponds to the trivial explicit system $(\RR^{m}_{\infty},\tau_{m})$, which justifies its name.


\section{Lie-B\"acklund equivalence for implicit systems}
Let us slightly adapt the notion of Lie-B\"acklund equivalence\footnote{In the terminology of \cite{FLMR-ieee}, different names have been introduced, which, in the author's opinion, are poorly matched: \emph{differential equivalence} (corresponding to \emph{endogenous transformations} and \emph{$\Phi$-related Cartan fields}) and Lie-B\"acklund equivalence (including time scalings into the previous endogenous transformations by replacing Cartan fields by Cartan distributions). In order to stress that they are two faces of the same coin, we propose to use \emph{Lie-B\"acklund equivalence} (resp.  \emph{orbital Lie-B\"acklund equivalence}) in place of \emph{differential equivalence} (resp. \emph{Lie-B\"acklund equivalence}) and \emph{Lie-B\"acklund isomorphisms} (resp.  \emph{orbital Lie-B\"acklund isomorphisms}) in place of  \emph{endogenous transformations} (resp. \emph{Lie-B\"acklund isomorphisms}).} of \cite{FLMR-ieee} in our implicit control system context:

Let us consider two regular implicit control systems $(\XX, \tau_{\XX},F)$, with $\XX= X\times \RR^{n}_{\infty}$, $\dim X=n$, $\tau_{\XX}$ the associated trivial Cartan field, and $\rk{\frac{\partial F}{\partial\dot{x}}}=n-m$, and $(\YY, \tau_{\YY},G)$, with $\YY= Y\times \RR^{p}_{\infty}$, $\dim Y=p$, $\tau_{\YY}$ the associated trivial Cartan field, and $\rk{\frac{\partial G}{\partial\dot{y}}}=p-q$.

Set $\XX_{0}=\{ \ol{x}\in \XX \vert L_{\tau_{\XX}}^{k}F(\ol{x})=0,~\forall k\geq 0 \}$ and $\YY_{0}=\{ \ol{y}\in \YY \vert L_{\tau_{\YY}}^{k}G(\ol{y})=0,~\forall k\geq 0 \}$. They are endowed with the topologies and differentiable structures induced by $\XX$ and $\YY$ respectively.

\begin{defn}\label{LBdef}
We say that two regular implicit control systems $(\XX, \tau_{\XX},F)$ and $(\YY, \tau_{\YY},G)$ are \emph{Lie-B\"acklund equivalent}  (or shortly L-B equivalent) at $(\ol{x}_{0},\ol{y}_{0})\in \XX_{0}\times \YY_{0}$ if, and only if:
\begin{itemize}
\item[(i)] there exist neighborhoods $\XXX_{0}$ of\,  $\ol{x}_{0}$ in $\XX_{0}$, and $\YYY_{0}$ of\,  $\ol{y}_{0}$ in $\YY_{0}$, and a one-to-one mapping $\Phi=(\varphi^{0},\varphi^{1},\ldots )\in C^{\infty}(\YYY_{0} ; \XXX_{0})$, with $C^{\infty}(\YYY_{0} ; \XXX_{0})$ inverse $\Psi$, satisfying $\Phi(\ol{y}_{0})=\ol{x}_{0}$ and such that the restrictions of the trivial Cartan fields ${\tau_{\YY}}_{\big| {\YYY_{0}}}$ and ${\tau_{\XX}}_{\big| {\XXX_{0}}}$are $\Phi$-related, namely $\Phi_{\ast}{\tau_{\YY}}_{\big| {\YYY_{0}}}={\tau_{\XX}}_{\big| {\XXX_{0}}}$;
\item[(ii)] the $C^{\infty}(\YYY_{0} ; \XXX_{0})$ inverse mapping $\Psi=(\psi^{0},\psi^{1},\ldots )$ is such that $\Psi(\ol{x}_{0})=\ol{y}_{0}$ and $\Psi_{\ast}{\tau_{\XX}}_{\big| {\XXX_{0}}}={\tau_{\YY}}_{\big| {\YYY_{0}}}$.
\end{itemize}
The mappings $\Phi$ and $\Psi$ are called \emph{mutually inverse Lie-B\"acklund isomorphisms} at $(\ol{x}_{0},\ol{y}_{0})$.

The two systems $(\XX, \tau_{\XX},F)$ and $(\YY, \tau_{\YY},G)$ are called \emph{locally L-B equivalent} if they are L-B equivalent at every pair $(\ol{x}, \Psi(\ol{x}))=(\Phi(\ol{y}),\ol{y})$ of an open dense subset $\ZZZ$ of $\XX_{0}\times \YY_{0}$, with $\Phi$ and $\Psi$ mutually inverse Lie-B\"acklund isomorphisms on $\ZZZ$.
\end{defn}
The next Proposition shows the equivalence of this definition to the one of \cite{FLMR-ieee} in the explicit context.
\begin{pr}
Given two systems $(\XX,\tau_{\XX},F)$ and $(\YY,\tau_{\YY},G)$ and two vector fields $f$ and $g$ compatible with $(\XX,\tau_{\XX},F)$ and $(\YY,\tau_{\YY},G)$ respectively. The corresponding explicit systems $\dot{x}=f(x,u)$ and $\dot{y}=g(y,v)$ are \emph{differentially equivalent} in the sense of \cite{FLMR-ieee} (or L-B equivalent as proposed in footnote \footnotemark[2]) at the pair $((x_{0},u_{0},\dot{u}_{0},\ldots),(y_{0},v_{0},\dot{v}_{0},\ldots))$, with $u_{0}$ such that  $\dot{x}_{0}=f(x_{0},u_{0})$ and $v_{0}$ such that $\dot{y}_{0}=g(y_{0},v_{0})$,  if, and only if, $(\XX,\tau_{\XX},F)$ and $(\YY,\tau_{\YY},G)$ are L-B equivalent according to Definition~\ref{LBdef}, at $(\ol{x}_{0},\ol{y}_{0})$ with $\ol{x}_{0}=(x_{0},\dot{x}_{0},\ldots)$ and $\ol{y}_{0}=(y_{0},\dot{y}_{0},\ldots)$.
\end{pr}
\begin{proof}
Let $\Phi$ and $\Psi$ satisfy \textit{(i)} and \textit{(ii)}. Since $\ol{g}$ is compatible with $(\YY, \tau_{\YY},G)$ and for \mbox{$\ol{y}\in \YYY_{0}$}, using the construction of Remark~\ref{implrem}, we have  $G(y,g(y,v))=0$ for all $v$ in a suitable open subset of $\RR^{q}$. Since, by assumption, $\ol{x}=\Phi(\ol{y})\in \XXX_{0} \subset \XX_{0}$, with $\Phi=(\varphi^{0},\varphi^{1},\ldots)$, we have \mbox{$x=\varphi^{0}(y,g(y,v),\frac{dg}{dt}(y,v,\dot{v}),\ldots)\triangleq \widetilde{\varphi}^{0}(y,\ol{v})$}, and $\ol{x}$ satisfies $L_{\tau_{\XX}}^{k}F(\ol{x})=0$ for all $k\geq 0$. Let $\ol{v}_{0}$ be such that $\dot{y}_{0}=g(y_{0},v_{0})$. Thus
\mbox{$\ol{x}_{0}=\Phi(y_{0},g(y_{0},v_{0}),\frac{d}{dt}g(y_{0},v_{0},\dot{v}_{0}),\ldots)$}.

Now, since $\ol{f}$ is compatible with $(\XX,\tau_{\XX},F)$ one has $u=\mu(x,\dot{x})$ again with the notation of Remark~\ref{implrem}, or  \mbox{$u=\mu(\widetilde{\varphi}^{0}(y,\ol{v}),\frac{d}{dt}\widetilde{\varphi}^{0}(y,\ol{v})) \triangleq \widetilde{\varphi}^{1}(y,\ol{v})$}. Using $\Phi_{\ast}{\tau_{\YY}}_{\big| {\YYY_{0}}}={\tau_{\XX}}_{\big| {\XXX_{0}}}$ yields 
$$\begin{array}{l}
\ds f(x,u)_{\big | (\widetilde{\varphi}^{0}(y,\ol{v}),\widetilde{\varphi}^{1}(y,\ol{v}))}=
L_{\tau_{\XX}}x_{\big | (\widetilde{\varphi}^{0}(y,\ol{v}),\widetilde{\varphi}^{1}(y,\ol{v}))}
=L_{\tau_{\YY}}\varphi^{0}(y, L_{\ol{g}}y, L_{\ol{g}}^{2}y,\ldots)\\
\ds = \sum_{j\geq 0}\sum_{i=1}^{p} \frac{\partial \varphi^{0}}{\partial y_{i}^{(j)}}L_{\ol{g}}\left( L_{\ol{g}}^{j}y_{i}\right)
= \sum_{j\geq 0}\sum_{i,k=1}^{p}g_{k} \frac{\partial \varphi^{0}}{\partial y_{i}^{(j)}}\frac{\partial L_{\ol{g}}^{j}y_{i}}{\partial y_{k}} + 
\sum_{j,l\geq 0} \sum_{k=1}^{q}\sum_{i=1}^{p}v_{k}^{(l+1)} \frac{\partial \varphi^{0}}{\partial y_{i}^{(j)}}\frac{\partial L_{\ol{g}}^{j}y_{i}}{\partial v_{k}^{(l)}}\\
\ds = \sum_{k=1}^{p}g_{k}\frac{\partial \widetilde{\varphi}_{0}}{\partial y_{k}} + \sum_{l\geq 0} \sum_{k=1}^{q}v_{k}^{(l+1)} \frac{\partial \widetilde{\varphi}^{0}}{\partial v_{k}^{(l)}} = (\widetilde{\varphi}^{0})_{\ast}\ol{g}(y,\ol{v}).
\end{array}$$  
Analogously, $\dot{u}= \frac{d}{dt}u=L_{\tau_{\XX}}\mu(\ol{x})= \frac{d}{dt}\widetilde{\varphi}^{1}(y,\ol{v})= \sum_{i=1}^{p}g_{i}(y,v)\frac{\partial \widetilde{\varphi}^{1}}{\partial y_{i}} + \sum_{j\geq 0} \sum_{i=1}^{q}v_{i}^{(j+1)} \frac{\partial \widetilde{\varphi}^{1}}{\partial v_{i}^{(j)}} = (\widetilde{\varphi}^{1})_{\ast}\ol{g}$, which proves that $(x,\ol{u})=\widetilde{\Phi}(y,\ol{v})$, with $\widetilde{\Phi}=(\widetilde{\varphi}^{0}, \widetilde{\varphi}^{1},\ldots)$, and $\ol{f}=\widetilde{\Phi}_{\ast}\ol{g}$, defining $\ol{u}_{0}$ by $(x_{0}, \ol{u}_{0})=\widetilde{\Phi}(y_{0},\ol{v}_{0})$. 

Symmetrically, we have $(y,\ol{v})=\widetilde{\Psi}(x,\ol{u})$ with $\widetilde{\Psi}=(\widetilde{\psi}^{0}, \widetilde{\psi}^{1},\ldots)$ and $\widetilde{\psi}^{0}(x,\ol{u})=\psi^{0}(x,f(x,u),\frac{df}{dt}(x,u,\dot{u}),\ldots)$, $\widetilde{\psi}^{1}(x,\ol{u})= \nu(\widetilde{\psi}^{0}(x,\ol{u}),\frac{d}{dt}\widetilde{\psi}^{0}(x,\ol{u}))$, $v=\nu(x,\dot{x})$, and thus $\widetilde{\Phi}$'s inverse is $\widetilde{\Psi}$ with $(y_{0},\ol{v}_{0})=\widetilde{\Psi}(x_{0},\ol{u}_{0})$, and $\ol{g}=\widetilde{\Psi}_{\ast}\ol{f}$, which proves that the explicit systems $(X\times\RR^{m}_{\infty},\ol{f})$ and $(Y\times\RR^{q}_{\infty},\ol{g})$ are locally L-B equivalent at $(x_{0}, \ol{u}_{0})$, $(y_{0},\ol{v}_{0})$ for  $\ol{u}_{0}$ and $\ol{v}_{0}$ suitably chosen, their choice depending on $f$ and $g$. The proof of the converse follows the same lines.
\end{proof}

An easy consequence of this definition is that L-B equivalence preserves equilibrium points, namely points $\ol{\widetilde{y}}=(\widetilde{y},0,0,\ldots)$ (resp. $\ol{\widetilde{x}}=(\widetilde{x},0,0,\ldots)$) such that $G(\widetilde{y},0)=0$ (resp. $F(\widetilde{x},0)=0$).

The following result is easily adapted from \cite{FLMR-ieee}:
\begin{pr}
If two regular implicit control systems $(\XX,\tau_{\XX},F)$ and $(\YY,\tau_{\YY}, G)$ are locally L-B equivalent, they have the same coranks, namely $m=q$.
\end{pr}


\section{Flatness and variational properties}
First recall from \cite{FLMR-ieee} that a system in explicit form is flat if, and only if, it is \mbox{L-B equivalent} to a trivial system. The reader may easily check that this definition is just a concise and precise restatement of the definition given in the introduction. The adaptation of this definition in our context is obvious: 

\begin{defn}\label{flatdef}
The implicit system $(\XX,\tau_{\XX},F)$ is \emph{flat} at $\ol{x}_{0}$ if, and only if, there exists \mbox{$\ol{y}_{0} \in \RR^{m}_{\infty}$} such that $(\XX,\tau_{\XX},F)$ is L-B equivalent, at $(\ol{x}_{0},\ol{y}_{0}) \in \XX_{0}\times \RR^{m}_{\infty}$, to the $m$-dimensional trivial implicit system $(\RR^{m}_{\infty},\tau_{m},0)$. In this case, the mutually inverse L-B isomorphisms $\Phi$ and $\Psi$ are called \emph{inverse trivializations}, (or \emph{uniformizations}, in reference to Hilbert's 22nd problem).

The system $(\XX,\tau_{\XX},F)$ is locally \emph{flat} if, and only if, there exists an open dense subset $\XXX_{0}$ of $\XX_{0}$ such that it is flat at every $\ol{x}_{0} \in \XXX_{0}$. 
\end{defn}

Otherwise stated, (\ref{implsys}) is flat at $\ol{x}_{0}$ if the local integral curves $t\mapsto \ol{x}(t)$ of (\ref{implsys}) around $\ol{x}_{0}$ are images by a smooth one-to-one mapping  $\Phi$, satisfying $\ol{x}_{0}=\Phi(\ol{y}_{0})$, of arbitrary curves in the coordinates $\ol{y}=(y_{1},\ldots,y_{m},\dot{y}_{1},\ldots,\dot{y}_{m},\ldots, y_{1}^{(k)},\ldots,y_{m}^{(k)},\ldots,)$ around  $\ol{y}_{0}$. In other words, for every curve $t\mapsto \ol{y}(t)$ in a suitable time interval $\II$, \mbox{$\ol{x}(t)=(x(t),\dot{x}(t),\ddot{x}(t),\ldots)=\Phi(\ol{y}(t))= (\varphi^{0}(\ol{y}(t)),\varphi^{1}(\ol{y}(t)),\varphi^{2}(\ol{y}(t)),\ldots)$} belongs to $\XX_{0}$ for all $t\in\II$ and thus $F(x(t),\dot{x}(t))=0$. Conversely, if $t\mapsto x(t)$ is an integral curve of $F(x,\dot{x})=0$, there exists a curve $t\mapsto y(t)$ in $C^{\infty}(\II;\RR^{m})$ such that $\ol{y}(t)=(y(t),\dot{y}(t),\ldots)=\Psi(\ol{x}(t))=(\psi^{0}(\ol{x}(t)),\psi^{1}(\ol{x}(t)),\ldots)$ for all $t\in \II$, namely $L_{\tau_{m}}\ol{y}(t)=L_{\tau_{\XX}}\Psi(\ol{x}(t))$ for all $t\in \II$.
Recall that $\Phi$ (resp. $\Psi$) depends only on a finite number of derivatives of $y$ (resp. $x$).

The extension of this remark to local flatness is straightforward.

Trivializations may be characterized in terms of the differential of $F$ as follows. 

A basis of the tangent space $\T_{\ol{x}}\XX$ of $\XX$ at a point $\ol{x}\in \XX$ consisting of the set of vectors \mbox{$\{ \frac{\partial}{\partial x_{i}^{(j)}} \vert i=1,\ldots,n, j\geq 0\}$}, a basis of the cotangent space $\T_{\ol{x}}^{\ast}\XX$ at $\ol{x}$ is therefore given by $\{ dx_{i}^{(j)} \vert i=1,\ldots,n, j\geq 0\}$ with \mbox{$<dx_{i}^{(j)},\frac{\partial}{\partial x_{k}^{(l)}}>=\delta_{i,k}\delta_{j,l}$}, $\delta_{i,k}$ being the Kronecker symbol (i.e. $\delta_{i,j}=1$ if $i=j$ and $\delta_{i,j}=0$ otherwise).
The differential of $F_{i}$ is  thus given by 
\begin{equation}\label{dFdG}
dF_{i}= \sum_{j=1}^{n} \left( \frac{\partial F_{i}}{\partial x_{j}}dx_{j} + \frac{\partial F_{i}}{\partial \dot{x}_{j}}d\dot{x}_{j}\right), \quad i=1,\ldots, n-m.
\end{equation}

Since smooth functions depend on a finite number of variables, their differential contains only a finite number of non zero terms. Accordingly, we define a 1-form on $\XXX$, an open dense subset of $\XX$, as a \emph{finite} linear combination of the  $dx_{i}^{(j)}$, with coefficients in $C^{\infty}(\XXX;\RR)$, or equivalently as a local $C^{\infty}$ section of $\T^{\ast}(X\times \RR^{n})$ (see e.g. \cite{Z}). The set of 1-forms on $\XXX$ is denoted by $\Lambda^{1}(\XXX)$. Clearly, the $dF_{i}$'s are elements of $\Lambda^{1}(\XX)$. 

Note that the shift property of $\ddt=L_{\tau_{\XX}}$ on coordinates extends to differentials by defining $\ddt dx_{i}$ as: $\ddt dx_{i}= d\dot{x}_{i}=d\ddt x_{i}$.
More generally, we define the Lie-derivative of a 1-form $\omega$ along a vector field $\ol{v}$ on $\XX$, which is a 1-form on $\XX$, denoted by $L_{\ol{v}}\omega$, as in the finite dimensional case by the Leibnitz rule:
$$\left < L_{\ol{v}}\omega,\ol{w} \right> = L_{\ol{v}}\left<\omega,\ol{w}\right> - \left< \omega,[\ol{v},\ol{w}]\right>$$ 
for every vector field $\ol{w}$ on $\XX$, where $[\ol{v},\ol{w}]$ is the Lie-bracket of $\ol{v}$ and $\ol{w}$.
Clearly, if $\omega=dx_{i}$, $\ol{v}=\tau_{\XX}$ and $\ol{w}=\frac{\partial}{\partial \dot{x}_{i}}$, we recover the previous formula, namely $L_{\tau_{\XX}}dx_{i}=d\dot{x}_{i}$.

If $\Phi$ is a smooth mapping  from $\YY$ to $\XX$, the definition of the image by $\Phi$ of a 1-form is the same as in the finite dimensional context: if $\omega\in \Lambda^{1}(\XX)$, $\omega(\ol{x})= \sum_{\scriptscriptstyle{\begin{array}{c} j\geq 0\\\mbox{\it{finite}}\end{array}}}\sum_{i=1}^{n} \omega_{j}^{i}(\ol{x})dx_{i}^{(j)}$, 
its (backward) image $\Phi^{\ast}\omega$ is the 1-form on $\YY$ defined by 
\begin{equation}\label{Phistarom}
\Phi^{\ast}\omega(\ol{y})=
\sum_{k,l}\sum_{i,j} \omega_{j}^{i}(\Phi(\ol{y}))\frac{\partial \varphi_{i}^{j}}{\partial y_{k}^{(l)}}(\ol{y})dy_{k}^{(l)}
\end{equation} 
where $\varphi_{i}^{j}$ is the $i,j$-th component of $\Phi$, namely $x_{i}^{(j)}=\varphi_{i}^{j}(\ol{y})$.

Note again that, since the functions $\varphi_{i}^{j}$ depend on a finite number of variables, the 1-form $\Phi^{\ast}\omega$ contains only a finite number of non zero terms and, according to $x_{i}^{(j)}=\frac{d^{j}x_{i}}{dt^{j}}$ and (\ref{Phistarom}), we have
\begin{equation}\label{dxij-eq}
\begin{array}{c}
\ds dx_{i}^{(j)}= \sum_{k,l}\frac{\partial \varphi_{i}^{j}}{\partial y_{k}^{(l)}} dy_{k}^{(l)} = L_{\tau_{\XX}}^{j}dx_{i} =
L_{\tau_{\YY}}^{j}\left(  \sum_{k,l}\frac{\partial \varphi_{i}^{0}}{\partial y_{k}^{(l)}}dy_{k}^{(l)} \right) \\
\ds =   \sum_{k,l}\sum_{r=0}^{j} {j\choose r}
\left( L_{\tau_{\YY}}^{r}\left( \frac{\partial \varphi_{i}^{0}}{\partial y_{k}^{(l)}} \right)\right) dy_{k}^{(l+j-r)}.
\end{array}
\end{equation}
where ${j\choose r} = \frac{j!}{r!(j-r)!}$ stands for the Bernoulli binomial coefficient.

\begin{thm}\label{cotanflat} The system $(\XX,\tau_{\XX},F)$ is flat at $\ol{x}_{0}$, with $\ol{x}_{0}\in \XX_{0}$, if, and only if, there exists $\ol{y}_{0}\in \RR^{m}_{\infty}$ and a local Lie-Bäcklund isomorphism $\Phi$ from $\RR^{m}_{\infty}$ to $\XX_{0}$ satisfying $\Phi(\ol{y}_{0})=\ol{x}_{0}$ and such that
\begin{equation}\label{cotanflat-eq}
\Phi^{\ast}dF=0.
\end{equation}
\end{thm}
\begin{proof}
Necessity: If the system $(\XX,\tau_{\XX},F)$ is flat at $\ol{x}_{0}$, we have \mbox{$F(\varphi_{0}(\ol{y}),\varphi_{1}(\ol{y}))=0$} for all $\ol{y}$ in a neighborhood of $\ol{y}_{0}$ in $\RR^{m}_{\infty}$ such that $\Phi(\ol{y}_{0})=\ol{x}_{0}$. For every $\lambda$ in a given interval $[0,\lambda_{0}[$ of $\RR$ and a sufficiently small time interval  $\II$ containing $0$, consider a smooth trajectory $t\mapsto y_{\lambda}(t)$ in a bounded neighborhood of $y_{0}$ in $\RR^{m}$, such that \mbox{$\sup \{ \Vert y^{(j)}_{\lambda}(t)\Vert \vert j\geq 0, t\in \II , \lambda \in [0,\lambda_{0}[ \}$} is finite, where $\Vert\cdot\Vert$ denotes the Euclidean norm of $\RR^{m}$, and set $\frac{\partial y_{\lambda}}{\partial \lambda}(t)_{\big| \lambda = 0} =\zeta(t)$, which exists on $\II$ by the Ascoli-Arzelà Theorem (see e.g. \cite{Y}). Next, we set $\ol{x}_{\lambda}(t)=\Phi(\ol{y}_{\lambda}(t))$ for $t\in \II$ and $\lambda \in [0, \lambda_{0}[$. We indeed have \mbox{$F(\varphi^{0}(\ol{y}_{\lambda}(t)),\varphi^{1}(\ol{y}_{\lambda}(t)))=0$} for all $t$ and $\lambda$ in their respective intervals and thus, differentiating with respect to $\lambda$ at $\lambda =0$, we get $\frac{\partial F}{\partial x}\frac{\partial \varphi^{0}}{\partial \ol{y}}(\ol{y}_{0}(t))\ol{\zeta}(t)+
\frac{\partial F}{\partial \dot{x}}\frac{\partial \varphi^{1}}{\partial \ol{y}}(\ol{y}_{0}(t))\ol{\zeta}(t)=0$ in $\II$.
At time $t=0$, noting $y_{\lambda}(0)=y_{\lambda,0}$, we have  $\frac{\partial F}{\partial x}\frac{\partial \varphi^{0}}{\partial \ol{y}}(\ol{y}_{0,0})\ol{\zeta}(0)+
\frac{\partial F}{\partial \dot{x}}\frac{\partial \varphi^{1}}{\partial \ol{y}}(\ol{y}_{0,0})\ol{\zeta}(0)=0$, and this expression is valid for every $\ol{y}_{0,0}$ in a neighborhood of $\ol{y}_{0}$ and $\ol{\zeta}(0)\in \T_{\ol{y}_{0,0}}\RR^{m}_{\infty}$. We have thus proved that the 1-form
 $\frac{\partial F}{\partial x}\frac{\partial \varphi^{0}}{\partial \ol{y}}d\ol{y}+
\frac{\partial F}{\partial \dot{x}}\frac{\partial \varphi^{1}}{\partial \ol{y}}d\ol{y}=\Phi^{\ast}dF$ vanishes on $\T_{\ol{y}_{0,0}}\RR^{m}_{\infty}$ for every $\ol{y}_{0,0}$ in a neighborhood of $y_{0}$, and therefore is identically zero in a neighborhood of $y_{0}$, which proves (\ref{cotanflat-eq}).

Sufficiency: assuming that there exists a locally smooth invertible mapping \mbox{$\Phi=(\varphi^{0},\varphi^{1},\ldots)\in C^{\infty}(\RR^{m}_{\infty};\XX)$} satisfying (\ref{cotanflat-eq}) with  $\Phi(\ol{y}_{0})=\ol{x}_{0}$, the 1-forms $\Phi^{\ast}dF_{i}$, $i=1,\ldots, n-m$, are obviously closed since they are the differentials of the functions $F_{i}\circ\Phi$, $i=1,\ldots, n-m$. Thus (\ref{cotanflat-eq}) implies that $F_{i}(\varphi^{0}(\ol{y}),\varphi^{1}(\ol{y}))=c_{i}$, $i=1,\ldots, n-m$,  with $c_{i}$ arbitrary constants. But since $\ol{x}_{0}\in \XX_{0}$ and $\ol{x}_{0}=\Phi(\ol{y}_{0})$, we have $F(x_{0},\dot{x}_{0})=0$ 
and $c_{i}=F_{i}(\varphi^{0}(\ol{y}_{0}),\varphi^{1}(\ol{y}_{0}))=F_{i}(x_{0},\dot{x}_{0})=0$, for $i=1,\ldots,n-m$.
Then, setting $\ol{x}=\Phi(\ol{y})=(\varphi^{0}(\ol{y}),\varphi^{1}(\ol{y}),\ldots)$, we get that $x=\varphi^{0}(\ol{y})$ (which depends on $y$ and a finite number of its derivatives) satisfies $F_{i}(x,\dot{x})=0$, $i=1,\ldots, n-m$ and that $\dot{x}=L_{\tau_{\XX}}x=L_{\tau_{\YY}}\varphi^{0}(\ol{y})=\varphi^{1}(\ol{y})$. Following the same lines for all derivatives of $x$, we have proved that $\Phi_{\ast}\tau_{\YY}=\tau_{\XX}$. Finally, since $\Phi$ is invertible with $C^{\infty}$ inverse $\Psi=(\psi^{0},\psi^{1},\ldots)$, it is immediately seen that $\ol{y}=\Psi(\ol{x})$ (and therefore $y=\psi^{0}(\ol{x})$ only depends on a finite number of derivatives of $x$) and that $\Psi_{\ast}\tau_{\XX}=\tau_{\YY}$, which proves the sufficiency and the proof is complete.
\end{proof}


\section{Flatness necessary and sufficient conditions}
\subsection{Preliminaries on polynomial matrices}
We now analyze condition (\ref{cotanflat-eq}) in greater details with the (mild) restriction that $F$ is \mbox{\emph{meromorphic}} on $\T X$. This restriction is motivated by the use of algebraic properties of polynomial matrices and of modules over a principal ideal ring of polynomials, this ring being itself formed over the field of meromorphic functions, as will be made clear immediately.

We also restrict the inverse trivializations $\Phi$ and $\Psi$ of definition~\ref{flatdef} to the class of meromorphic functions.

In matrix notations and using indifferently $\ddt$ for $L_{\tau_{\XX}}$ or $L_{\tau_{\YY}}$ (the context being unambiguous), according to (\ref{dxij-eq}), we have: 
$$\begin{array}{l}
\ds \Phi^{\ast} dF = \sum_{j\geq 0}\sum_{i=1}^{m}
\left(  \frac{\partial F}{\partial x}\frac{\partial \varphi^{0}}{\partial y_{i}^{(j)}} + \frac{\partial F}{\partial \dot{x}}\frac{\partial \varphi^{1}}{\partial y_{i}^{(j)}} 
 \right)  dy_{i}^{(j)} \\
 \hspace{1.1cm} \ds =  \sum_{j\geq 0}\sum_{i=1}^{m}
 \left( \frac{\partial F}{\partial x} \frac{\partial \varphi^{0}}{\partial y_{i}^{(j)}} dy_{i}^{(j)} + 
\frac{\partial F}{\partial \dot{x}} \frac{d}{dt}
\left( \frac{\partial \varphi^{0}}{\partial y_{i}^{(j)}} dy_{i}^{(j)}  \right) \right)\\
\hspace{1.1cm} \ds =
\left( \frac{\partial F}{\partial x} + \frac{\partial F}{\partial \dot{x}} \frac{d}{dt} \right)_{\big| \ol{x}=\Phi(\ol{y})}
\left( \sum_{i=1}^{m}\sum_{j\geq 0} 
\frac{\partial \varphi^{0}}{\partial y_{i}^{(j)}} \frac{d^{j}}{dt^{j}} dy_{i} \right).
\end{array}
$$
Since $\varphi^0$ depends only on a finite number of derivatives of $y$, there exists a finite integer $j^{\ast}$ such that there exists $i\in \{1,\ldots,m\}$ for which $\frac{\partial \varphi^{0}}{\partial y_{i}^{(j^{\ast})}} \neq 0$ and 
$\frac{\partial \varphi^{0}}{\partial y_{i}^{(j)}} = 0$ for every $i\in \{1,\ldots,m\}$ and every $j \geq j^{\ast}$. This integer $j^{\ast}$ indeed represents the maximum degree in $\ddt$ of the polynomials
$\sum_{j\geq 0} \frac{\partial \varphi^{0}}{\partial y_{i}^{(j)}} \frac{d^{j}}{dt^{j}}$, for $i=1,\ldots m$, and is denoted by $j^{\ast} =\ord{\varphi^{0}}$. Therefore
$$
\sum_{j\geq 0} \frac{\partial \varphi^{0}}{\partial y_{i}^{(j)}} \frac{d^{j}}{dt^{j}} = \sum_{j = 0}^{\ord{\varphi^{0}}} \frac{\partial \varphi^{0}}{\partial y_{i}^{(j)}} \frac{d^{j}}{dt^{j}}, \quad i=1,\ldots m.
$$

Introducing the following polynomial matrices where the indeterminate is the differential operator $\ddt$:
\begin{equation}\label{PF-PPhi}
P(F) = \frac{\partial F}{\partial x} + \frac{\partial F}{\partial \dot{x}} \frac{d}{dt} ,\quad
P(\varphi^{0}) = \sum_{j = 0}^{\ord{\varphi^{0}}} \frac{\partial \varphi^{0}}{\partial y^{(j)}} \frac{d^{j}}{dt^{j}}
\end{equation}
with $P(F)$ (resp. $P(\varphi^{0})$) of size $(n-m)\times n$ (resp. $n\times m$), (\ref{cotanflat-eq}) reads:
\begin{equation}\label{PFpoly}
\Phi^{\ast} dF_{\big| \ol{y}}= P(F)_{\big| \Phi(\ol{y})} P(\varphi^{0})_{\big| \ol{y}} \ dy= 0.
\end{equation}
or equivalently
\begin{equation}\label{PFpoly1}
\Phi^{\ast} dF_{\big| \Psi(\ol{x})}= P(F)_{\big| \ol{x}}  P(\varphi^{0})_{\big| \Psi(\ol{x})} dy= 0.
\end{equation}
Clearly, the entries of these matrices are polynomials of the differential operator $\ddt$ whose coefficients are meromorphic functions from $\XX$ to $\RR$. We denote by $\kk$ the field of meromorphic functions from  $\XX$ to $\RR$ and by $\kk[\ddt]$ the principal ideal ring of polynomials of $\ddt=L_{\tau_{\XXX}}$ with coefficients in $\kk$.

We may also consider the field of meromorphic functions from $\YY=\RR^{m}_{\infty}$ to $\RR$. In this case, the notations $\kk_{\YY}$ and $\kk_{\YY}[\ddt]$, with $\ddt =L_{\tau_{\YY}}$, will replace the previous ones.

Recall that $\kk[\ddt]$ is non commutative, even if $n=1$, as shown by the following example:  denoting, with abusive notations, by $x$ the 0th order operator $a\mapsto xa$ for every $a\in \kk$, we have, for $a\neq 0$, $\left( \ddt\cdot x -x\cdot \ddt\right)(a)=\dot{x}a+x\dot{a}-x\dot{a}=\dot{x}a \neq 0$, or $ \ddt\cdot x -x\cdot \ddt=\dot{x}$. 

For  arbitrary integers $p$ and $q$, let us denote by $\MM_{p,q}[\ddt]$ the module of $p\times q$ matrices over $\kk[\ddt]$ (see e.g. \cite{Co,KS} for a detailed presentation of modules over non commutative rings).
Recall that, since in general the inverse of a polynomial is not a polynomial, for arbitrary $p\in \NN$, the inverse of a square invertible matrix of $\MM_{p,p}[\ddt]$ doesn't generally belong to $\MM_{p,p}[\ddt]$. Matrices whose inverse belong to $\MM_{p,p}[\ddt]$ are called \emph{unimodular matrices} and their set is denoted by $\UU_{p}[\ddt]$. It forms a normal\footnote{$A^{-1} \UU_{p}[\ddt] A = \UU_{p}[\ddt]$ for every invertible $A\in \MM_{p,p}[\ddt]$} subgroup of the group generated by invertible matrices of $\MM_{p,p}[\ddt]$.

Recall from \cite[Chap.8]{Co} (see the Annex~\ref{polymat-sec}) the following fundamental result:
\begin{thm}[Smith decomposition (or diagonal reduction)]\label{polymat-th}
Given a matrix $M\in \MM_{p,q}[\ddt]$, there exist matrices
$V\in \UU_{p}[\ddt]$ and
$U\in \UU_{q}[\ddt]$
such that
\begin{equation}\label{M-decomp}
VMU=\left( \Delta_{p}, 0_{p,q-p}\right)~\mbox{\textrm{if~}} p\leq q, \quad 
VMU=\left( \begin{array}{c}\Delta_{q}\\ 0_{p-q,q}\end{array}\right)~\mbox{\textrm{if~}} p\geq q
\end{equation}
where
$0_{p,q-p}$ (resp. $0_{p-q,q}$) is the $q\times (q-p)$ (resp. $(p-q)\times q$) matrix whose entries are all zeros, and with
$\Delta_{p}$  a $p\times p$ (resp. $\Delta_{q}$ a $q\times q$) diagonal matrix whose
diagonal elements, $(\delta_{1},\ldots,\delta_{\sigma},0,\ldots,0)$, are such
that
$\delta_{i}$ is a non zero $\ddt$-polynomial for $i=1,\ldots,\sigma$, and is a  divisor
of~$\delta_{j}$ for all~$\sigma\geq j\geq i$. Moreover, $\Delta_{p}$ (resp. $\Delta_{q}$ is unique up to multiplication by a regular diagonal matrix in $\kk^{p\times p}$ (resp. $\kk^{q\times q}$).
\end{thm}
The above unimodular matrices $U$ and $V$ are indeed non unique. We say that $U\in \rsm{M}$ (resp. $V\in \lsm{M}$) if there exists $V_U$ (resp. $U_V$) such that the pair $(U,V_U)$ (resp. $(U_V,V)$) satisfies (\ref{M-decomp}).

\bigskip

Here, $P(F) \in \MM_{n-m,n}[\ddt]$. According to Theorem~\ref{polymat-th}, it admits the \emph{Smith decomposition}
\begin{equation}\label{PFdecomp}
V   P(F)   U =\left( \Delta , 0_{n-m,m} \right)
\end{equation}
with $V \in \UU_{n-m}[\ddt]$, $U \in \UU_{n}[\ddt]$ and $\Delta \in \MM_{n-m,n-m}[\ddt]$.

\begin{defn}
Given a matrix $M\in \MM_{p,q}[\ddt]$, we say that $M$ is \emph{hyper-regular} if, and only if, its Smith decomposition gives $\left( \Delta_p , 0_{p,q-p} \right) = \left(I_{p},0_{p,q-p} \right)$ if $p\leq q$ and
$\left( \begin{array}{c}\Delta_{q}\\ 0_{p-q,q}\end{array}\right) = \left(\begin{array}{c} I_{q}\\0_{p-q,q}\end{array} \right)$ if $p\geq q$. 
\end{defn}
Note that a square matrix $M\in \MM_{p,p}[\ddt]$ is hyper-regular if, and only if, it is unimodular, i.e. $M\in \UU_{p}[\ddt]$.

\bigskip

Consider a point $\ol{x}\in \XX$ and its projection $x$ on the original manifold $X$. Following \cite{Fl-scl},   the \emph{variational module} ${\cal M}$ of (\ref{implsys}) at $\ol{x}$ is the finitely generated module constructed as follows:
consider $\xi=(\xi_{1},\ldots,\xi_{n})$ a non zero but otherwise arbitrary element of the tangent space $\T_{x}X$ and denote by $[\xi]$ the $\kk[\ddt]$-module generated by the components of $\xi$. We also denote by $[P(F)\xi]$ the submodule of $[\xi]$ generated by the components of $P(F)\xi$. Then, the 
 \emph{variational module} ${\cal M}$ of (\ref{implsys}) at $\ol{x}$ is the quotient module $[\xi]/[P(F)\xi]$.
According to \cite{Fl-scl} (see also e.g. \cite{Co}), ${\cal M}$ can be decomposed into the following direct sum: 
$${\cal M}={\cal T} \oplus {\cal F}$$ 
where the uniquely defined module ${\cal T}$ is torsion and $\cal F$ is free. $\cal F$ is unique up to isomorphism.\\

It is readily seen that:
\begin{itemize}
\item ${\cal T}$ is the $\kk[\ddt]$-module generated by the components of  $U \zeta$, in ${\cal M}$, with $\zeta$ satisfying \mbox{$\delta_{i,i}\zeta_{i}=0$} for some $i= 1, \ldots, n-m$, and $\zeta_{n-m+1}= \cdots = \zeta_{n}=0$, for $U\in \rsm{P(F)}$,
\item and ${\cal F}$ is a free $\kk[\ddt]$-module generated by the components of $U \left(\begin{array}{cc} 0_{n-m, n-m}& 0_{n-m, m}\\ 0_{m, n-m}& I_{m}\end{array}\right) \zeta$ for every $\zeta$ whose components are in ${\cal M}$, and $U \in \rsm{P(F)}$.
\end{itemize}

For the trivial system $\YY=\RR^{m}_{\infty}$, its variational module ${\cal M}_{\YY}$ at an arbitrary point $\ol{y}\in \RR^{m}_{\infty}$ is identified to the tangent space $\T_{y}\RR^{m}$.

Clearly, if $P(F)$ is hyper-regular, since the polynomial degree of every diagonal element of $\Delta = I_{m}$ is 0, we immediately deduce that ${\cal T}=\{ 0\}$ and that ${\cal M}$ is free. 

The next subsection establishes a link between flatness of the system corresponding to (\ref{implsys}), F-controllability of its variational module, and hyper-regularity of $P(F)$.

\subsection{Flatness and controllability}
Recall from  \cite{Fl-scl} that the variational system of (\ref{implsys}) at $\ol{x}$ is said to be \emph{F-controllable}\footnote{The prefix F- (which may be equally understood as Free or Fliess) has been added to avoid confusions with other linear or nonlinear controllability concepts. Comparisons with such notions are not addressed in this paper.} if, and only if, its associated module is free.
 
\begin{pr}\label{controllability-prop}
If system $(\XX,\tau_{\XX},F)$ is locally flat at $\ol{x}_{0}$, its variational system at every point of a neighborhood of $\ol{x}_{0}$ is F-controllable and $P(F)$ is hyper-regular in this neighborhood.
\end{pr}
\begin{proof}
Let $(\ol{x},\ol{y})$ belong to a neighborhood of $(\ol{x}_{0},\ol{y}_{0})$ in $\XX_{0}\times \RR^{m}_{\infty}$ with $\ol{x}=\Phi(\ol{y})$, or equivalently $\ol{y}=\Psi(\ol{x})$, $\Phi$ and $\Psi$ being mutually inverse trivializations.  

Assume that the variational module ${\cal M}$ associated to (\ref{implsys}) is torsion, i.e. ${\cal T} \neq \{0\}$. Thus according to what precedes, there exists a non zero vector  $\zeta \in \T_{\ol{x}}\XX$ whose components are  torsion elements,  and a pair of unimodular matrices $U$ and $V$ such that $P(F)U\zeta = V^{-1}\Delta \zeta = 0$.  Denote by $\theta = \Psi_{\ast}\zeta = \frac{\partial \Psi}{\partial \ol{y}} \zeta$. Clearly, $\theta \in \T_{\ol{y}}\RR^{m}_{\infty}$ and its components belong to ${\cal M}_{\YY}$, the free $\kk_{\YY}[\ddt]$-module associated to the trivial system. We also have $\zeta = \Phi_{\ast} \theta = \frac{\partial \Phi}{\partial \ol{x}} \theta$, and thus $V^{-1}\Delta \frac{\partial \Phi}{\partial \ol{x}} \theta = 0$. Since $V$ and $\frac{\partial \Phi}{\partial \ol{x}}$ are locally invertible, the polynomial matrix $V^{-1}\Delta \frac{\partial \Phi}{\partial \ol{x}}$ doesn't identically vanish, which proves that $\theta$ is a torsion element of ${\cal M}_{\YY}$, which contradicts the fact that ${\cal M}_{\YY}$ is free. Therefore, local flatness implies F-controllability.

To prove that F-controllability implies hyper-regularity of $P(F)$, assume that $P(F)$ is not hyper-regular. Without loss of generality, according to the fact that the diagonal elements $d_{i,i}$ of $\Delta$ divide $d_{j,j}$ for $1\leq i \leq j \leq n-m$, we may assume that all the diagonal elements $d_{i,i}$ are in $\kk$ but the last one, $d_{n-m,n-m}$, a polynomial of degree greater than or equal to 1 in $\ddt$. Consider the homogeneous differential equation $d_{n-m,n-m}\zeta_{n-m} =0$. Since its coefficients are meromorphic, it has a unique local non zero meromorphic solution for every non zero initial condition and therefore the components of the non zero $n$-dimensional vector $\xi=U\zeta$ with $\zeta = (0, \ldots, 0, \zeta_{n-m}, 0, \ldots, 0)^{T}$ and $\zeta_{n-m}$ at the $(n-m)$th place, belong to ${\cal M}$ since $P(F)\xi= P(F)U\zeta= \Delta \zeta =0$, and contain at least one torsion element since $\Delta U^{-1}\xi =\Delta \zeta = 0$, which contradicts the  F-controllability. The proof is complete.
\end{proof}



\subsection{Algebraic characterization of the differential of a trivialization}
From now on, we assume that $P(F)$ is hyper-regular in a neighborhood $\XXX_{0}$ of $\ol{x}_{0}\in \XX_{0}$, since otherwise the system cannot be flat.

Thus, there exist $V\in \lsm{P(F)}$ and $U\in \rsm{P(F)}$ such that 
\begin{equation}\label{PFdecompI}
V  P(F)  U =\left(I_{n-m}, 0_{n-m,m} \right).
\end{equation}

In this framework, the set of all polynomial matrices $P(\varphi^{0})_{\big| \Psi(\ol{x})}\in \MM_{n,m}[\ddt]$ satisfying (\ref{PFpoly}), or (\ref{PFpoly1}), can be completely characterized. We first solve the matrix equation:
\begin{equation}\label{PFpolyform}
P(F) \Theta = 0
\end{equation}
where the entries of $\Theta\in  \MM_{n,m}[\ddt]$ are not supposed to be gradients of some function $\varphi^{0}$. Next, we show that to every such $\Theta$ is associated a flat output $\omega$ of the variational system, defined by the relation $dx= \Theta \omega$, which can be inverted: $\omega= \Xi dx$ for some polynomial matrix $\Xi$. We will then study the integrability aspects of the latter equation, or more precisely characterize the existence of an integrating factor $M$ such that $dy=M\omega$, and of $\Psi$ such that $y=\Psi(\ol{x})$ with $d\Psi =M\Xi dx$, in order to obtain a flat output $y$, if it exists, of the nonlinear system $(X\times\RR^{n}_{\infty},\tau_{X},F)$ in section~\ref{integ-subsec}.
\begin{lem}\label{Smith-lem}
Let $U\in \rsm{P(F)}$ and denote by
\begin{equation}\label{Uhatdef}
\widehat{U} =U \left( \begin{array}{c} 0_{n-m,m}\\I_{m} \end{array} \right).
\end{equation}
Every hyper-regular matrix $\Theta\in  \MM_{n,m}[\ddt]$ satisfying (\ref{PFpolyform}) is given by
\begin{equation}\label{PPhidecomp}
\Theta = \widehat{U} W
\end{equation}
with $W \in \UU_{m}[\ddt]$ arbitrary. 
\end{lem}
\begin{proof}
First, note that the set of hyper-regular matrices $\Theta\in  \MM_{n,m}[\ddt]$ satisfying (\ref{PFpolyform}) is non empty. 
Using Theorem~\ref{polymat-th} of Appendix~\ref{Smith-app}, in a suitable neighborhood where $\rk{\frac{\partial F}{\partial \dot{x}}}=n-m$, relation (\ref{PFdecompI}) implies that
$$VP(F)U \left( \begin{array}{c} 0_{n-m,m}\\I_{m} \end{array} \right)= VP(F)\widehat{U} =
\left(I_{n-m}, 0_{n-m,m} \right)
\left( \begin{array}{c} 0_{n-m,m}\\I_{m} \end{array} \right) = 0.$$
Thus, $P(F)\widehat{U}=0$, which means that $\widehat{U}$ is solution of (\ref{PFpolyform}).

Let $\Theta$ be an arbitrary hyper-regular solution of equation (\ref{PFpolyform}) and $U$ a unimodular matrix in $\rsm{P(F)}$. Using again (\ref{PFdecompI}), we get
\begin{equation}\label{PFpolyform1}
V  P(F)  U    U^{-1}    \Theta  =\left( I_{n-m} , 0_{n-m,m}\right)  U^{-1}   \Theta  =0
\end{equation} 
from which we deduce that
$U^{-1}  \Theta  =\left( \begin{array}{c} 0_{n-m,m}\\I_{m} \end{array} \right)   W$, where $W$ is an arbitrary $m\times m$ polynomial matrix, hence (\ref{PPhidecomp}), using (\ref{Uhatdef}). The hyper-regularity of $\Theta$ and the invertibility of $U$  immediately imply the hyper-regularity of $W$. But  since $W$ is a square ($m\times m$) hyper-regular matrix, it is unimodular, hence the lemma. 
\end{proof}
\begin{rem}
Solutions $\Theta$ of (\ref{PFpolyform}) are fully characterized by the fact that their range is equal to the kernel of $P(F)$, i.e. we have the (local) short exact sequence of modules
$$
0 \longrightarrow \left(\Lambda^{1}(\XX)\right)^{m}  \stackrel{\Theta}{\longrightarrow} \left(\Lambda^{1} (\XX)\right)^{n} \stackrel{P(F)}{\longrightarrow}  \left(\Lambda^{1} (\XX)\right)^{n-m}  \longrightarrow  0 
$$
This s is in accordance with the fact that, as a consequence of this Lemma, a solution $\Theta$ does not depend on a particular choice of $U\in \rsm{P(F)}$: according to (\ref{PPhidecomp}), for two different choices $U_1$ and $U_2$, there corresponds two different choices $W_1$ and $W_2$ such that $\Theta= \widehat{U}_1 W_1 = \widehat{U}_2 W_2$. 
\end{rem}
\begin{lem}[``Inversion'' of $\Theta$]\label{Smith2-lem}
For every $U\in \rsm{P(F)}$ and $Q \in \lsm{\widehat{U}}$, there exists $Z\in \UU_{m}[\ddt]$ such that
\begin{equation}\label{Uhatdecomp}
Q \Theta =\left( \begin{array}{c} I_{m}  \\ 
0_{n-m,m}\end{array}\right) Z.
\end{equation}
Moreover, splitting $Q$ into the two blocks:
\begin{equation}\label{Qtilde-hat}\widetilde{Q}= \left(\begin{array}{cc}I_{m},&0_{m,n-m}\end{array}\right)Q, \quad \widehat{Q}=\left(0_{n-m,m},I_{n-m}\right)Q
\end{equation}
we have 
\begin{equation}\label{Uhatdecomp-Z}
\widetilde{Q}\Theta =Z, \quad \widehat{Q}\Theta=0_{n-m,m}
\end{equation}
and $\widehat{Q}$ is equivalent to $P(F)$, i.e. there exists a matrix $L \in \UU_{n-m}[\ddt]$ such that $P(F)=L\widehat{Q}$.

In particular, there exists a $Q$ corresponding to $\Theta= \widehat{U}$ and $Z=I_{m}$, denoted by $Q_{0}$:
\begin{equation}\label{Q0def-Id}
Q_{0}= \left(\begin{array}{cc}0_{m,n-m}&I_{m}\\I_{n-m}&0_{n-m,m}\end{array}\right)
U^{-1} \in \lsm{\widehat{U}},
\end{equation}
and the corresponding blocks
\begin{equation}\label{Q0tilde}\widetilde{Q}_{0}= \left(\begin{array}{cc}I_{m},&0_{m,n-m}\end{array}\right)Q_{0}, \quad \widehat{Q}_{0}=\left(0_{n-m,m},I_{n-m}\right) Q_{0}
\end{equation}
satisfy
\begin{equation}\label{Uhatdecomp-Id}
Q_{0}\widehat{U} =\left( \begin{array}{c} I_{m}\\0_{n-m,m}\end{array}\right), \quad \widetilde{Q}_{0}\widehat{U} =I_{m}, \quad \widehat{Q}_{0}\widehat{U}=0_{n-m,m}
\end{equation}
with $\widehat{Q}_{0}$ equivalent to $P(F)$, i.e. there exists $L_{0} \in \UU_{n-m}[\ddt]$ such that $P(F)=L_{0}\widehat{Q}_{0}$.
\end{lem}
\begin{proof}
Since $\widehat{U}$, by (\ref{Uhatdef}), is hyper-regular,
there exist $Q\in \UU_{n}[\ddt]$ and $R\in \UU_{m}[\ddt]$ such that $Q\widehat{U}R=\left(\begin{array}{c}I_{m}\\0_{n-m,m}\end{array}\right)$. Using (\ref{PPhidecomp}), (\ref{Uhatdecomp}) follows by setting $R^{-1}W=Z\in \UU_{m}[\ddt]$, and we immediately deduce (\ref{Uhatdecomp-Z}) by left multiplying (\ref{Uhatdecomp}) by $\left(\begin{array}{cc}I_{m},&0_{m,n-m}\end{array}\right)$ and $\left(0_{n-m,m},I_{n-m}\right)$ successively.
Moreover, comparing $\widehat{Q}\Theta=0_{n-m,m}$, in (\ref{Uhatdecomp-Z}), with $P(F)\Theta=0_{n-m,m}$, and taking account of the hyper-regularity of $\Theta$, $\widehat{Q}$ and $P(F)$, the existence of a unimodular matrix $L$ such that $P(F)=L\widehat{Q}$ is proven. 

In (\ref{Q0def-Id}), the proof that $Q_{0}\in \lsm{\widehat{U}}$ is immediate by direct computation of the product $Q_{0}\widehat{U}$, which also proves the left equality of (\ref{Uhatdecomp-Id}). The remaining equalities of (\ref{Uhatdecomp-Id}) also follow by direct computation and the last assertion of the Lemma is just a restatement of the previous one in the particular case $\Theta= \widehat{U}$ and $Z=I_{m}$, which completes the proof of the lemma. 
\end{proof}
\begin{rem}
The Smith decomposition used in Lemmas~\ref{Smith-lem} and \ref{Smith2-lem}, whose algorithm is given in Annex~\ref{polymat-sec}, may be easily implemented in a computer algebra system in this non commutative context. However, as far as the number of operations is concerned, this algorithm might be improvable. But since we only focus attention on the existence of the matrices $\widehat{U}$ and $Q$,  this problem is not addressed here. 
\end{rem}

Consider now $U\in \rsm{P(F)}$, $Q_{0}$ defined by (\ref{Q0def-Id}), and $\widetilde{Q}_{0}$ given by (\ref{Q0tilde}).
Let us denote by $Q_{i}^{j}=\sum_{\alpha \geq 0}Q_{i,\alpha}^{j}\frac{d^{\alpha}}{dt^{\alpha}}$ the $(i,j)$-th polynomial entry of $Q_{0}\in \lsm{\widehat{U}}$. We define the $m$-dimensional vector 1-form $\omega$ by
\begin{equation}\label{omeg-def}
\omega(\ol{x})= \left( \begin{array}{c} \omega_{1}(\ol{x})\\\vdots\\\omega_{m}(\ol{x}) \end{array}\right) \triangleq \widetilde{Q}_{0}dx_{\big| \XXX_{0}} 
= 
\left( \begin{array}{c} {\sum_{j=1}^{n}\sum_{\alpha \geq 0}Q_{1,\alpha}^{j}(\ol{x}) dx_{j}^{(\alpha)}}_{\big| \XXX_{0}} \\
\vdots\\
{\sum_{j=1}^{n}\sum_{\alpha \geq 0}Q_{m,\alpha}^{j}(\ol{x})dx_{j}^{(\alpha)}}_{\big| \XXX_{0}} 
\end{array}\right) 
\end{equation}
the restriction to the neighborhood $\XXX_{0}$ of $\XX_{0}$ meaning that the expressions of the right-hand side of (\ref{omeg-def}) are evaluated for every $\ol{x}\in \XXX_{0} \subset \XX_{0}$ satisfying $L_{\tau_{\XX}}^{\alpha}F=0$ and every $dx_{j}^{(\alpha)}$ such that $dL_{\tau_{\XX}}^{\alpha}F=0$ for all $\alpha$.

Since $\widetilde{Q}_{0}$ is hyper-regular, the forms $\omega_{1},\ldots,\omega_{m}$ are independent.

Let us recall that, adapting the general definition~\ref{flatdef} to the linear context, in the language of modules, the components of a flat output of the variational system are simply the elements of a basis of the free module ${\cal M}$ (see \cite{FLMR-ijc,Chyzak_05}).

The two previous lemmas may be reformulated as:
\begin{cor}[Variational flat outputs]\label{flatvar-cor}
The vector 1-form $\omega=\widetilde{Q}_{0}dx$, defined by (\ref{omeg-def}) with $Q_{0}$ given by (\ref{Q0def-Id}) and $\widetilde{Q}_{0}$ by (\ref{Q0tilde}), is a flat output of the variational system
and we have $dx= \widehat{U}\omega$, or
\begin{equation}\label{dxflat}
dx_{i}= \sum_{j=1}^{m}\sum_{\alpha=0}^{\ord{U}}\widehat{U}_{i,\alpha}^{j} \ \omega_{j}^{(\alpha)},\quad i=1,\ldots,n.
\end{equation}
Similarly, denoting by $\ord{\widetilde{Q}_{0}}$ the polynomial order of $\widetilde{Q}_{0}$ with respect to $\ddt$ and $\#(\widetilde{Q}_{0})$ the maximal number of derivatives of $x$ which the entries of $\widetilde{Q}_{0}$ depend on, we have $\ord{\widetilde{Q}_{0}}\leq (n-1)\ord{U}$, $\#(\widetilde{Q}_{0})\leq \#(U) +(n-1)\ord{U}$ and there exist integers $\ord{\Gamma}$ and $\#(\Gamma)$ satisfying 
$$\ord{\Gamma}\leq \#(U) + n(\ord{U}), \quad 
\#(\Gamma) \leq 2\#(U) + (n-1)\ord{U}$$
 and meromorphic functions $\Gamma_{i,\alpha,\beta}^{j,k}$, $j,k=1,\ldots, m$, $\alpha, \beta =0,\ldots, \ord{\Gamma}$, depending at most of $(x,\dot{x},\ldots, x^{(\#(\Gamma))})$, such that
\begin{equation}\label{structureconst-eq}
d\omega_{i} = \sum_{\alpha,\beta = 0}^{\ord{\Gamma}} \sum_{j,k=1}^{m} \Gamma_{i,\alpha,\beta}^{j,k} \ \omega_{j}^{(\alpha)}\wedge \omega_{k}^{(\beta)}.
\end{equation}
Finally, the set of flat outputs of the variational system is equal to $\UU_{m}[\ddt] \omega$.
\end{cor}
\begin{proof}
Assume that $dx$ satisfies $P(F)dx = 0$. By Lemma~\ref{Smith2-lem}, since $\widetilde{Q}_{0}\widehat{U}=I_{m}$, we have $\widetilde{Q}_{0}\left( dx- \widehat{U}\omega\right)=0$, which proves that there exists a 1-form $\zeta \in \ker{\widetilde{Q}_{0}}$ such that $dx=\widehat{U}\omega + \zeta$.
But, again by Lemma~\ref{Smith2-lem}, we have $P(F)=L_{0}\widehat{Q}_{0}$ with $L_{0}\in \UU_{n-m}[\ddt]$ and $\widehat{Q}_{0}=\left(0_{n-m,m},I_{n-m}\right) Q_{0}$, which yields $0=L_{0}^{-1}P(F)dx = \widehat{Q}_{0}dx = \widehat{Q}_{0}\widehat{U}\omega +  \widehat{Q}_{0}\zeta$. Since $\widehat{Q}_{0}\widehat{U}=0$, we immediately get that $\zeta \in \ker{\widetilde{Q}_{0}} \cap \ker{\widehat{Q}_{0}}= \{ 0 \}$ and $dx=\widehat{U}\omega$.
Therefore, $dx$ can be expressed as a function of $\omega$ and its successive derivatives, i.e. $dx= \sum_{\alpha=0}^{\ord{U}}\widehat{U}_{\alpha}\omega^{(\alpha)}$, the finiteness of the integer $\ord{U}$ resulting from the Smith decomposition algorithm, and conversely, by (\ref{omeg-def}), $\omega$ may be expressed as a function of $dx$ and its successive derivatives, which proves that $\omega$ is a flat output of the variational system.

Moreover, since $\widetilde{Q}_{0}= \left(\begin{array}{cc}I_{m}&0_{m,n-m}\end{array}\right)
\left(\begin{array}{cc} 0_{n-m,m}&I_{m}\\I_{n-m}&0_{n-m,m}\end{array}\right) U^{-1}$. and since the degree of the inverse $U^{-1}$ of $U\in \UU_{n}[\ddt]$ is bounded by $(n-1)\ord{U}$ (see \cite{Ritt,Oll-aecc,KMP}), we conclude that  $\ord{\widetilde{Q}_{0}}\leq (n-1)\ord{U}$.

It is left as an exercise that, as an immediate consequence of the previous inequality, we have 
$\#(U^{-1}) \leq \#(U) +(n-1)\ord{U}$, and therefore $\#(\widetilde{Q}_{0}) \leq \#(U) +(n-1)\ord{U}$.

Taking the exterior derivative of both sides of (\ref{omeg-def}), the $i$th component of the vector 2-form $d\omega$ is given by
\begin{equation}\label{d-omega-i}
d\omega_{i}= \sum_{\gamma=0}^{\ord{\widetilde{Q}_{0}}} \sum_{\eta=0}^{\#(\widetilde{Q}_{0})} \sum_{l,r=1}^{n}
\sum_{j,k=1}^{m}
\frac{\partial \widetilde{Q}_{i,\gamma}^{l}}{\partial x_{r}^{(\eta)}} \;\; dx_{r}^{(\eta)} \wedge dx_{l}^{(\gamma)}.
\end{equation}
Differentiating (\ref{dxflat}), we get
\begin{equation}\label{dxr-eta}
dx_{r}^{(\eta)}=  \sum_{j=1}^{m} \sum_{\rho=0}^{\ord{U}} \sum_{\alpha=\rho}^{\rho + \eta} 
{\eta\choose \alpha - \rho}
\left( \widehat{U}_{r,\rho}^{j} \right)^{(\eta + \rho - \alpha)} \omega_{j}^{(\alpha)}
\end{equation}
where ${\eta\choose \alpha - \rho}$ is the binomial coefficient $\frac{\eta !}{(\alpha - \rho)! (\eta - \alpha + \rho)!}$.
Thus, the combination of (\ref{dxr-eta}) and (\ref{d-omega-i}), after reordering the summations, yields (\ref{structureconst-eq}), with:
\begin{equation}\label{d-omega-i-Gamma}
\begin{aligned}
\Gamma_{i,\alpha,\beta}^{j,k}= 
\sum_{l,r=1}^{n} 
\sum_{\rho=0}^{\min(\ord{U},\alpha)}  \sum_{\sigma=0}^{\min(\ord{U},\beta)} 
\sum_{\eta=\max(0,\alpha-\rho)}^{\#(\widetilde{Q}_{0})} 
\sum_{\gamma=\max(0,\beta-\sigma)}^{\ord{\widetilde{Q}_{0}}} 
&{\eta\choose \alpha - \rho}
{\gamma\choose \beta - \sigma}\cdot
\\
&\cdot \frac{\partial \widetilde{Q}_{i,\gamma}^{l}}{\partial x_{r}^{(\eta)}}
\left( \widehat{U}_{r,\rho}^{j} \right)^{(\eta + \rho - \alpha)} \left( \widehat{U}_{l,\sigma}^{k} \right)^{(\gamma + \sigma - \beta)}
\end{aligned}
\end{equation}
for all $\alpha = 0,\ldots,\ord{U}+\#(\widetilde{Q}_{0})$ and $\beta= 0,\ldots,\ord{U}+\ord{\widetilde{Q}_{0}}$, which proves (\ref{structureconst-eq}) with $\ord{\Gamma}$ satisfying 
\begin{equation}\label{o(Gamma)}
\ord{\Gamma} \leq \max \left(\ord{U}+\#(\widetilde{Q}_{0}), \ord{U}+\ord{\widetilde{Q}_{0}} \right) = \#(U) + n(\ord{U}). 
\end{equation}
Moreover, as an immediate consequence of (\ref{d-omega-i-Gamma}), the maximal number $\#(\Gamma)$ of derivatives of $x$ in the 
$\{\Gamma_{i,\alpha,\beta}^{j,k}\}$'s satisfies the inequality
\begin{equation}\label{[Gamma]}
\#(\Gamma) \leq \max \left(\#(U)+ \#(\widetilde{Q}), \#(U) + \ord{\widetilde{Q}} \right)= 2\#(U) + (n-1)\ord{U},
\end{equation} 
hence the result.

Finally, any other flat output $\kappa$ of the variational system, i.e. whose components form a basis of the free module ${\cal M}$, may be deduced from $\omega$ by  $\kappa= M\omega$ for some $M\in \UU_{m}[\ddt]$ and conversely, given an arbitrary $M\in \UU_{m}[\ddt]$, it is readily seen that $\kappa=M\omega$ is also a basis of ${\cal M}$. Therefore, the set of flat outputs of the variational system is $\UU_{m}[\ddt] \omega$, which completes the proof.
\end{proof}
\begin{rem} The functions $\Gamma_{i,\alpha,\beta}^{j,k}$ of (\ref{structureconst-eq}) are non unique as will be seen later in Proposition~\ref{sol-mu}.
\end{rem}


\subsection{Integrability}\label{integ-subsec}
Let us also recall that, if $\tau_{1},\ldots,\tau_{r}$ are given independent 1-forms in $\Lambda^{1}(\XXX_{0})$, the $\kk[\ddt]$-ideal $\TTT$ generated by $\tau_{1},\ldots,\tau_{r}$, for an arbitrary integer $r$, is the set of all combinations with coefficients in  $\kk[\ddt]$ of forms $\eta\wedge \tau_{i}$ with $\eta$ an arbitrary form on $\XXX_{0}$ of arbitrary degree and $i=1,\ldots,r$.

Note that if another set of independent 1-forms $\kappa_{1},\ldots,\kappa_{s}$ of $\TTT$ is a generator of  $\TTT$, then $s=r$ and there exists a unimodular matrix $H\in \UU_{r}[\ddt]$ such that $\tau=H\kappa$, where $\tau=\left(\tau_{1},\ldots,\tau_{r}\right)^{T}$ (the superscript $^{T}$ means transposition) and $\kappa=\left(\kappa_{1},\ldots,\kappa_{r}\right)^{T}$.
\begin{defn}
We say that the $\kk[\ddt]$-ideal $\TTT$ generated by $\tau_{1},\ldots,\tau_{r}$ is \emph{strongly closed} if, and only if, there exists a matrix $M\in \UU_{m}[\ddt]$ such that $d(M\tau)=0$.
\end{defn}
This definition is indeed independent of the choice of generators since if $\kappa$ is another vector of generators of $\TTT$, we have $\tau=H\kappa$ for $H\in \UU_{r}[\ddt]$, and since $d(M\tau)=0$, we have $d(MH\kappa)=0=d(M'\kappa)$, with $M'=MH \in \UU_{r}[\ddt]$, which proves our assertion.

\begin{thm}\label{CNSthm}
A necessary and sufficient condition for system $(\XX,\tau_{\XX},F)$ to be flat at $(\ol{x}_{0},\ol{y}_{0})$ (over the class of meromorphic functions) is that
there exist $U\in \rsm{P(F)}$ and $Q\in\lsm{\widehat{U}}$, with $\widehat{U}$ given by (\ref{Uhatdef}),  such that the $\kk[\ddt]$-ideal $\Omega$ generated by the 1-forms $\omega_{1},\ldots,\omega_{m}$ defined by (\ref{omeg-def}) is strongly closed in $\XXX_{0}$.
\end{thm}
\begin{proof}
\textbf{Necessity}:  If system (\ref{implsys}) is flat at $(\ol{x}_{0},\ol{y}_{0})$, there exists $\Phi=\left( \varphi_{0},\varphi_{1},\ldots\right)$ meromorphic from a neighborhood $\YYY_{0}$ of $\ol{y}_{0}$ in $\RR^{m}_{\infty}$ to a neighborhood $\XXX_{0}$ of $\ol{x}_{0}$ in $\XX_{0}$ and one-to-one, such that $\ol{x}=\Phi(\ol{y})$ implies $F(\varphi_{0}(\ol{y}),\varphi_{1}(\ol{y}))=0$ and $L_{\tau_{\YY}}^{k}(F\circ \Phi)=0$ for all $k\geq 0$. According to Theorem~\ref{cotanflat}, $P(\varphi_{0})$ satisfies (\ref{PFpoly}), or (\ref{PFpoly1}), and is hyper-regular: proceeding as in Proposition~\ref{controllability-prop}, its kernel $P(\varphi_{0})^{-1}(\{0\})$ cannot contain torsion elements since otherwise it would contradict the freeness property of the module ${\cal M}$ (corresponding to the variational system of the trivial system $(\RR^{m}_{\infty},\tau_{m},0)$).  Thus, in virtue of Lemma~\ref{Smith-lem}, there exists $W\in \UU_{m}[\ddt]$ such that 
\begin{equation}\label{Pphi0nec}
P(\varphi_{0})_{\big| \Psi(\ol{x})}=\widehat{U}_{\big | \ol{x}} W_{\big | \ol{x}}.
\end{equation}
By Lemma~\ref{Smith2-lem},  we have $Q_{0}\widehat{U}= I_m$ and thus
\begin{equation}\label{XtilPphi0nec}
(\widetilde{Q}_{0})_{\big| \ol{x}}  P(\varphi_{0})_{\big| \Psi(\ol{x})}= W_{\big| \ol{x}} .
\end{equation}
Taking the exterior derivative of $x=\varphi_{0}(\ol{y})$ yields $dx=P(\varphi_{0}) dy$, and thus, according to (\ref{XtilPphi0nec}) and (\ref{omeg-def}),  $\omega (\Phi(\ol{y}))= (\widetilde{Q}_{0})_{\big| \Phi(\ol{y})} dx= (\widetilde{Q}_{0})_{\big| \Phi(\ol{y})}   P(\varphi_{0})_{\big| \ol{y}} dy= W_{\big| \Phi(\ol{y})} dy$, or 
\begin{equation}\label{Z-1Xdx}
\left(W_{\big| \ol{x}}\right)^{-1}  (\widetilde{Q}_{0})_{\big| \ol{x}}dx= \left(W_{\big| \ol{x}}\right)^{-1} \omega_{\big| \ol{x}} =dy.
\end{equation}
The forms $\omega_{1},\ldots,\omega_{m}$ are independent, generate the ideal $\Omega$, and setting $M=\left(W \right)^{-1}$ in (\ref{Z-1Xdx}) we get $M\omega=dy$. Thus,
taking the exterior derivative of both sides yields  $d(M\omega)=0$, which proves that $\Omega$ is strongly closed.

\textbf{Sufficiency}: Assume that there exist $U\in \rsm{P(F)}$ and $Q_{0}\in\lsm{\widehat{U}}$ given by (\ref{Q0def-Id}), such that $\Omega$, generated by the forms $\omega_{1},\ldots,\omega_{m}$ defined by  (\ref{omeg-def}), is strongly closed in a neighborhood $\XXX_{0}$ of $\ol{x}_{0}$ in $\XX_{0}$. Let $M\in\UU_{m}[\ddt]$ be such that $d(M\omega)=0$. Setting $\eta=M\omega$, the 1-forms $\eta_{1},\ldots,\eta_{m}$ finitely generate $\Omega$, are independent in $\XXX_{0}$, and, when expressed in the basis $dx_{i}^{(j)}$, contain only a finite number of terms whose coefficients depend on a finite number of derivatives of $x$. In the corresponding finite dimensional manifold, we have $d\eta_{i}=0$, $i=1,\ldots,m$, and, by Poincar\'e's Lemma, there locally exists a mapping $\psi_{0}\in C^{\infty}(\XXX_{0};\RR^{m})$ such that $d\psi_{0}=\eta=M\omega=M\widetilde{Q}_{0}dx$. In addition $\psi_{0}$ is a meromorphic function of its arguments since its differential is, according to the previous relation. 

Denoting by $y=\psi_{0}(\ol{x})$ for all $\ol{x}\in\XXX_{0}$ and $\Psi=\left( \psi_{0},\ddt\psi_{0}, \frac{d^{2}}{dt^{2}}\psi_{0},\ldots \right)=\left( \psi_{0},\psi_{1}, \psi_{2}, \ldots\right)$,
we have to prove that  $\Psi$ is a trivialization. 

Since $d\psi_{0}=\eta=M\omega=M\widetilde{Q}_{0}dx$, and since, according to Corollary~\ref{flatvar-cor}, $dx=\widehat{U}\omega$, we get $dx= \widehat{U}M^{-1}d\psi_{0}$.
If we denote by $\sigma_{i}$ the highest polynomial degree of the entries of the $ith$ column of $\widehat{U}M^{-1}$, we must have $n\leq m+ \sigma_{1}+\ldots +\sigma_{m}$ since otherwise this would contradict the surjectivity of $\widehat{U}M^{-1}$, considered as the matrix $\Xi$ whose entries are the $(\widehat{U}M^{-1})_{i,j}^{k}$'s, mapping an open subset of $\RR^{\sigma_{1}+1}\times\ldots\times \RR^{\sigma_{m}+1}$ to an open subset of $\RR^{n}$. In addition, if we note $\sigma=\max(\sigma_{i}\vert i=1,\ldots,m)$ and $\ol{y}^{\sigma}= \left( y_{1}^{(0)},\ldots,y_{1}^{(\sigma_{1})},\ldots,y_{m}^{(0)},\ldots,y_{m}^{(\sigma_{m})}\right)$,  $\Xi = \frac{\partial x}{\partial \ol{y}^{\sigma}}$ has rank $n$. 
Then the implicit system 
\begin{equation}\label{psi0implsys}
\begin{array}{l}
\ds y=\psi_{0}(\ol{x})\\
\ds \dot{y}=\psi_{1}(\ol{x})\\
\vdots\\
\ds y^{(\sigma)}=\psi_{\sigma}(\ol{x})
\end{array}
\end{equation}
has rank $n$ with respect to $x$ since its Jacobian matrix is a pseudo-inverse of $\Xi$. Hence, by the implicit function Theorem, a local solution to (\ref{psi0implsys}) is given by $x=\varphi_{0}( y,\ldots, y^{(\sigma)},\dot{x},\ldots,x^{(\rho)})$ for a suitable $\rho\in \NN$. But differentiating $\varphi_{0}$, using the fact that $dF=0$, or equivalently $P(F)dx=0$, and comparing with (\ref{PPhidecomp}), we find that $\varphi_{0}$ is independent of $(\dot{x},\ldots,x^{(\rho)})$, or $x=\varphi_{0}(\ol{y})$. It results that $\Phi=(\varphi_{0},\ddt\varphi_{0},\ldots)$ is the inverse trivialization of $\Psi$ which completes the proof of the Theorem.
\end{proof}


\subsection{Extension of the exterior derivative to polynomial differential forms}
Before giving a characterization of strongly closed ideals, we need to introduce some notations and tools. First, let us denote by $\Lambda^{p}(\XX)$ the module of all $p$-forms on $\XX$ . Note that the elements of $\Lambda^{1}(\XX)$ may be identified with $1$-forms on $X$ whose coefficients are in $\kk[\ddt]$ by the formula $\sum_{\alpha \geq 0}\sum_{i=1}^{n} \kappa_{\alpha}^{i}dx_{i}^{(\alpha)}= \sum_{i=1}^{n} \sum_{\alpha \geq 0} \left( \kappa_{\alpha}^{i} \frac{d^{\alpha}}{dt^{\alpha}} \right) dx_{i}$. We also denote by $\left( \Lambda^{p}(\XX)\right)^{m}$ the space of all $m$-dimensional vector $p$-forms on $\XX$, by $\left( \Lambda (\XX)\right)^{m}$ the space of all the $m$-dimensional vector forms of arbitrary degree on $\XX$, and by
$\lin{q}{ \left( \Lambda(\XX) \right)^{m}} \triangleq \lin{}{ \left( \Lambda^{p}(\XX) \right)^{m}, \left( \Lambda^{p+q}(\XX) \right)^{m} , p \geq 1 }$, the space of all linear operators from $\left( \Lambda^{p}(\XX) \right)^{m}$ to $\left( \Lambda^{p+q}(\XX)\right)^{m}$ for all $p\geq 1$, where $\lin{}{{\mathcal P},{\mathcal Q}}$ denotes the set of linear mappings from a given space ${\mathcal P}$ to a given space ${\mathcal Q}$.\\ 
Typically, an element $\mu \in \lin{q}{ \left( \Lambda(\XX) \right)^{m}}$ is an $m\times m$ matrix whose $(i,j)$-th entry reads
\begin{equation}\label{muform}
\mu_{i}^{j} = \sum_{\alpha\geq 0}\mu_{i,\alpha}^{j}\wedge\frac{d^{\alpha}}{dt^{\alpha}}, \quad i,j=1,2,\ldots,m
\end{equation}
where $\mu_{i,\alpha}^{j} \in \Lambda^{q}(\XX)$ for every $i,j=1,\ldots,m$, $\alpha\geq 0$. Therefore $\mu$ may be identified with a $\ddt$ polynomial whose coefficients are $m\times m$ matrices with entries in $\Lambda^{q}(\XX)$ ($q$-forms on $\XX$), i.e. $\mu \in \MM_{m,m}(\Lambda^{q}(\XX))[\ddt]$.
Thus, for every $\kappa \in \left( \Lambda^{p}(\XX) \right)^{m}$, the $i$th component of $\mu\kappa$ is given by
$$\left( \mu\kappa\right)_{i} = \sum_{\alpha \geq 0}\sum_{j=1}^{m} \mu_{i,\alpha}^{j} \wedge  L_{\tau_{\XX}}^{\alpha}\kappa_{j}$$ which is a $(p+q)$-form for every $i=1,\ldots,m$.

We define the operator $\dgot$, as an extension of the exterior derivative operator to unimodular matrices, by:
\begin{equation}\label{dgoth}
\dg{H}\kappa = d(H\kappa) -Hd\kappa
\end{equation}
for all  $m$-dimensional vector $p$-form $\kappa$ in $\left( \Lambda^{p}(\XX)\right)^{m}$, all $p\geq 1$, and all $H\in \UU_{m}[\ddt]$.
Note that (\ref{dgoth}) uniquely defines $\dg{H}$ as an element of $\lin{1}{ \left( \Lambda(\XX) \right)^{m}}$.

We can prolong $\dgot$ for all $\mu \in \lin{q}{ \left( \Lambda(\XX)\right)^{m}}$ and for all $\kappa \in \left( \Lambda^{p}(\XX)\right)^{m}$ and all $p\geq 1$ by the formula:
\begin{equation}\label{dgothprlg}
\dg{\mu}\kappa = d(\mu\  \kappa) - (-1)^{q}\mu\ d\kappa.
\end{equation}

For $\mu$ given by (\ref{muform}), it is straightforward to check that the $(i,j)$th entry $\dg{\mu}_{i}^{j}$ of $\dg{\mu}$ is given by
\begin{equation}\label{dgotmu}
\dg{\mu}_{i}^{j}= \sum_{\alpha \geq 0} d\mu_{i,\alpha}^{j} \wedge \frac{d^{\alpha}}{dt^{\alpha}}.
\end{equation}

By (\ref{dgothprlg}), we have, for every $i=1,\ldots, m$,
$$\begin{aligned}
\left( \dg{\mu}\kappa \right)_{i} =& d\left(\sum_{\alpha \geq 0}\sum_{j=1}^{m} \mu_{i,\alpha}^{j} \wedge  L_{\tau_{\XX}}^{\alpha}\kappa_{j}\right)- (-1)^{q}\left(  \sum_{\alpha \geq 0}\sum_{j=1}^{m} \mu_{i,\alpha}^{j} \wedge  L_{\tau_{\XX}}^{\alpha}d\kappa_{j}\right) \\
=& \left( \sum_{\alpha \geq 0}\sum_{j=1}^{m} d\mu_{i,\alpha}^{j} \wedge  L_{\tau_{\XX}}^{\alpha}\kappa_{j} + (-1)^{q}\sum_{\alpha \geq 0}\sum_{j=1}^{m} \mu_{i,\alpha}^{j} \wedge  L_{\tau_{\XX}}^{\alpha}d\kappa_{j} \right) - (-1)^{q}\left(  \sum_{\alpha \geq 0}\sum_{j=1}^{m} \mu_{i,\alpha}^{j} \wedge  L_{\tau_{\XX}}^{\alpha}d\kappa_{j}\right)
\end{aligned}
$$
hence the result.

The operator $\dgot$ enjoys the following properties:
\begin{pr}\label{complex-pr}
For all $\mu \in \lin{q}{ \left( \Lambda(\XX)\right)^{m}}$, all $q\in \NN$ and all $\kappa \in \left( \Lambda^{p}(\XX)\right)^{m}$, with $p\geq 1$ arbitrary, we have
\begin{equation}\label{complex}
\dg{\dg{\mu}}\kappa = 0.
\end{equation}
In other words $\dgot^{2}=0$, i.e. $\dgot$ is a complex.
\end{pr}
\begin{proof}
According to (\ref{dgothprlg}), replacing $\mu$ by $\dg{\mu}\in  \lin{q+1}{ \left( \Lambda(\XX)\right)^{m}}$, we get
\begin{equation}\label{d2goth}
\dg{\dg{\mu}}\kappa = d\left( \dg{\mu}\kappa\right) - (-1)^{q+1} \dg{\mu}d\kappa.
\end{equation}
Since, again with (\ref{dgothprlg}),  $d\left( \dg{\mu}\kappa\right) = d^{2}\left( \mu\ \kappa\right)
- (-1)^{q}d\left( \mu\ d\kappa\right)$ and, since $d^{2}=0$, we have
$$\begin{array}{l}
\ds \dg{\dg{\mu}}\kappa = - (-1)^{q}\left( \dg{\mu}\ d\kappa + (-1)^{q} \mu \ d^{2}\kappa \right) - (-1)^{q+1}\dg{\mu}\ d\kappa\\
\ds  = -(-1)^{q}\dg{\mu}\ d\kappa + (-1)^{q}\dg{\mu}\ d\kappa = 0
\end{array}$$ 
the result is proven.
\end{proof}

We also have:
\begin{pr}\label{dgmu}
For all $H\in \UU_{m}[\ddt]$ and all
$\mu \in \lin{q}{ \left( \Lambda(\XX)\right)^{m} }$, with arbitrary $q\geq 0$, the following relation holds:
\begin{equation}\label{dgothproduit}
\dg{H}\mu+H\dg{\mu}=\dg{H\mu}.
\end{equation}
In particular, if $\mu=-H^{-1}\dg{H}\in \lin{1}{ \left( \Lambda(\XX)\right)^{m} }$, we have
\begin{equation}\label{dgothmu}
\dg{\mu}=\mu^{2}.
\end{equation}
\end{pr}
\begin{proof} Let $H\in \UU_{m}[\ddt]$, $\kappa \in \left( \Lambda^{p}(\XX)\right)^{m}$ and 
$\mu \in \lin{q}{ \left( \Lambda(\XX)\right)^{m}}$,  with $q\geq 0$ arbitrary. (\ref{dgoth}) and (\ref{dgothprlg}), yield
$$\dg{H}\mu\kappa = d(H\mu\kappa)-Hd(\mu\kappa), \quad 
d(\mu\kappa)=\dg{\mu}\kappa + (-1)^{q}\mu d\kappa$$ 
or
$$\dg{H}\mu\kappa = d(H\mu\kappa)-H\dg{\mu}\kappa - (-1)^{q}H\mu d\kappa.$$
In other words:
$$\left( \dg{H}\mu + H\dg{\mu} \right) \kappa  = d(H\mu\kappa) - (-1)^{q}H\mu d\kappa = \dg{H\mu}\kappa.$$
This relation being valid for all $\kappa \in \left( \Lambda^{p}(\XX)\right)^{m}$ and all $p\geq 1$, we immediately deduce (\ref{dgothproduit}).

If now $\mu=-H^{-1}\dg{H}\in \lin{1}{ \left( \Lambda^{p}(\XX)\right)^{m}}$, we get $- H\mu = \dg{H}$ and thus, according to what precedes,
$$\dg{H\mu}=\dg{H}\mu+H\dg{\mu}=- \dgot^{2}(H)=0$$
or $H\dg{\mu} = - \dg{H}\mu$, or also $\dg{\mu} = - H^{-1}\dg{H}\mu= \mu^{2}$, which completes the proof.
\end{proof}

\bigskip

\subsection{A computable flatness characterization}

A characterization of the strong closedness condition on $\Omega$ is given by the next:
\begin{thm}\label{movingframe-thm}
The $\kk[\ddt]$-ideal $\Omega$ generated by the components of the vector 1-form $\omega$ defined by (\ref{omeg-def}) is strongly closed in $\XXX_{0}$ (or, equivalently, the system $(\XX,\tau_{\XX},F)$ is flat) if, and only if, there exists $\mu \in \lin{1}{ \left( \Lambda(\XX)\right)^{m} }$, and a matrix $M\in \UU_{m}[\ddt]$ such that
\begin{equation}\label{movingframe}
d\omega = \mu \ \omega, \qquad \dg{\mu} =  \mu^{2}, \qquad \dg{M}=- M\mu
\end{equation}
with the notation $\mu^2 =\mu\mu$.

In addition, if (\ref{movingframe}) holds true, a flat output $y$ is obtained by integration of $dy=M\omega$.
\end{thm}
\begin{proof}
If $d(M\omega)=0$, according to (\ref{dgoth}), we have
$\dg{M}\omega = -Md\omega$
or $d\omega = -M^{-1}\dg{M}\omega$. Setting $\mu= -M^{-1}\dg{M} \in \lin{1}{ \left( \Lambda(\XX)\right)^{m} }$, which is equivalent to the last relation of (\ref{movingframe}), we have
$d\omega = -M^{-1}\dg{M}\omega = \mu\  \omega$.
The conditions (\ref{movingframe}) are thus immediately deduced from (\ref{dgothmu}).

Conversely, from the last identity of (\ref{movingframe}), we get $-M^{-1}\dg{M}=\mu$. Its combination with the first one yields
$d\omega=-M^{-1}\dg{M}\omega$, or $Md\omega = -\dg{M}\ \omega$ and, according to (\ref{dgoth}), we immediately get that $d(M\omega) =0$, i.e. $\Omega$ is strongly closed.

Finally, if there exists a matrix $M\in \UU_{m}[\ddt]$ such that $d(M\omega)=0$, By Poincaré's Lemma, there exist $m$ functions $y_{1},\ldots,y_{m}$ such that $dy=M\omega$, which completes the proof.
\end{proof}
\begin{rem} Condition (\ref{movingframe}) may be seen as a generalization in the framework of manifolds of jets of infinite order of the well-known moving frame structure equations (see e.g. \cite[Chap 6, §3]{CCL}). Proceeding with this analogy, $\mu$ may be interpreted as a generalized curvature, and the fact that $\dgot(\mu)-\mu^{2}=0$ by the absence of (generalized) torsion.
\end{rem}
\begin{rem} The necessary and sufficient conditions (Theorem 4) of Chetverikov \cite{Cht-diffiety} (see also
 \cite{Cht}) show some similarities with (\ref{movingframe}) of our Theorem~\ref{movingframe-thm}, if we put aside the fact that the former results are obtained for explicit systems and in the $C^{\infty}$ context. More precisely, the two first conditions of (\ref{movingframe}) are similar to Chetverikov's conditions (A) and (B), where the operator $R$ is the analog of our $\mu$ (note that the basis $\omega$ and $R$ in \cite{Cht-diffiety} may depend on the $du_{i}^{(j)}$'s, the differentials of the control variables and their successive derivatives, which are eliminated in our context). However, condition (C) of this Theorem, which states that the operators generated by the iterated generalized symmetries applied to $R$ must have a bounded polynomial degree, is different from ours and seems to be difficult to verify in practice.
\end{rem}
\begin{pr}\label{sol-mu}
The general matrix solution $\mu=\left( \mu_{i}^{k}\right)_{i,k=1,\ldots,m}$ of $d\omega= \mu \omega$, with $\omega$ defined by (\ref{omeg-def}), is given by
\begin{equation}\label{mugen}
\begin{array}{c}
\ds \mu_{i}^{k} = \sum_{j=1}^{m} \sum_{\alpha,\beta=0}^{\ord{\mu}} \left( \Gamma_{i,\alpha,\beta}^{j,k} +\nu_{i,\alpha,\beta}^{j,k}\right) \omega_{j}^{(\alpha)} \wedge \frac{d^{\beta}}{dt^{\beta}}, \vspace{0.5em}\\
\mbox{\rm{with~}} \left\{ \begin{array}{l}
\ds \nu_{i,\alpha,\beta}^{j,k} = \nu_{i,\beta,\alpha}^{k,j} \quad \forall i,j,k=1,\ldots, m, \;\; \forall \alpha, \beta= 0,\ldots, \ord{\mu}, \; \alpha \neq \beta \; \mbox{\rm{or~}} j\neq k, \\
\ds \nu_{i,\alpha,\alpha}^{k,k}~\mbox{\rm{arbitrary},~} \quad \forall i,k=1,\ldots, m, \;\;  \forall \alpha = 0,\ldots, \ord{\mu}.
\end{array} \right.
\end{array}
\end{equation}
the integer $\ord{\mu}$ being arbitrary but otherwise finite and satisfying $\ord{\mu}\geq \ord{\Gamma}$, the $\Gamma_{i,\alpha,\beta}^{j,k}$'s being given by (\ref{structureconst-eq}), and where the $\nu_{i,\alpha,\beta}^{j,k}$'s are meromorphic functions depending at most on $\#(\mu)$ successive derivatives of $x$, with $\#(\mu)$ a finite integer such that $\#(\mu)\geq \#(\Gamma)$.
\end{pr}
\begin{proof}
Since $\mu$ is a matrix of the form $\mu=\sum_{\alpha=0}^{\ord{\mu}} \widetilde{\mu}_{\alpha}\wedge \frac{d^{\alpha}}{dt^{\alpha}}$
with $\widetilde{\mu}_{\alpha}\in \MM_{m,m}\left(\Lambda^1(\XX)\right)$ for every $\alpha$ and arbitrary finite integer $\ord{\mu}$, and since $\omega$ is a flat output of the variational system, choosing $\ord{\mu}$ large enough, every $\widetilde{\mu}_{\alpha}$ may be expressed in the coframe $\left \{ \omega, \dot{\omega}, \ldots, \omega^{(\ord{\mu})} \right \}$, yielding
$\widetilde{\mu}_{i,\alpha}^{j}= \sum_{k=1}^{m} \sum_{\beta=0}^{\ord{\mu}} \widetilde{\mu}_{i,\beta,\alpha}^{k,j} \ \omega_{k}^{(\beta)} \wedge \frac{d^{\alpha}}{dt^{\alpha}}$, where we have noted $\widetilde{\mu}_{i,\alpha}^{j}$ the $(i,j)$th entry of $\widetilde{\mu}_{\alpha}$. By linearity, every solution $\mu$ of $d\omega= \mu \omega$ is given by $\mu= \mu_0 + \mu_1$ with 
$\mu_0$ satisfying $\mu_{0} \omega = 0$ and $\mu_1=  \sum_{j=1}^{m} \sum_{\alpha,\beta=0}^{\ord{\mu}}  \Gamma_{i,\alpha,\beta}^{j,k}  \;  \omega_{j}^{(\alpha)} \wedge \frac{d^{\beta}}{dt^{\beta}}$, according to (\ref{structureconst-eq}), and 
$\mu_{0}$ is given by 
$\mu_{0} = \left( \nu_{i}^{j}\right)_{i,j=1,\ldots,m}$ with
$\nu_{i}^{j} =\sum_{k=1}^{m} \sum_{\alpha,\beta=0}^{\ord{\mu}} \nu_{i,\beta,\alpha}^{k,j} \ \omega_{k}^{(\beta)} \wedge \frac{d^{\alpha}}{dt^{\alpha}}$.
Moreover, since, for every $i=1,\ldots,m$: 
$$\begin{aligned}
 0 =  \left(\mu_{0} \omega\right)_{i} &= \sum_{j,k=1}^{m} \sum_{\alpha,\beta=0}^{\ord{\mu}} \nu_{i,\beta,\alpha}^{k,j} \ \omega_{k}^{(\beta)} \wedge \omega_{j}^{(\alpha)}\\
&= \sum_{0\leq \beta < \alpha \leq \ord{\mu}} \sum_{j,k=1}^{m}  \left( \nu_{i,\beta,\alpha}^{k,j} - \nu_{i,\alpha,\beta}^{j,k} \right) \omega_{k}^{(\beta)} \wedge \omega_{j}^{(\alpha)}  
+  \sum_{\beta=0}^{\ord{\mu}} \sum_{1\leq j < k \leq m}  \left( \nu_{i,\beta,\beta}^{k,j} - \nu_{i,\beta,\beta}^{j,k} \right) \omega_{k}^{(\beta)} \wedge \omega_{j}^{(\beta)} 
\end{aligned}
$$ 
and since the $\omega_{k}^{(\beta)} \wedge \omega_{j}^{(\alpha)}$'s form a basis of $\T^{\ast}(\XXX_{0})\wedge \T^{\ast}(\XXX_{0})$, we deduce that
$\mu_{i,\alpha,\beta}^{j,k} = \mu_{i,\beta,\alpha}^{k,j}$ for all $i,j,k = 1, \ldots, m$ and all
$\alpha, \beta= 0,\ldots, \ord{\mu}$, with $\alpha \neq \beta$ or $j\neq k$, hence the result.
\end{proof}
\begin{cor}\label{existmuM-cor}
The differential system (\ref{movingframe})
is algebraically closed.
\end{cor}
\begin{proof}
Applying the operator $\dgot$ to the two first equations of (\ref{movingframe}), we indeed obtain $\dgot(\mu\omega)=d(\mu\omega)=\dgot(\mu)\omega -\mu d\omega = \mu^{2}\omega-\mu^{2}\omega = 0$, and $\dgot(\mu^{2})=\dgot(\mu)\mu- \mu\dgot(\mu)=\mu^{3}-\mu^{3}=0$. Then, again applying $\dgot$ to $\dgot(M)=-M\mu$, we get $-\dgot(M)\mu-M\dgot(\mu)=M\mu^{2}-M\mu^{2}=0$ and the result is proven.
\end{proof}
\begin{rem}\label{existM}
Corollary~\ref{existmuM-cor} guarantees that no other exterior differential equations can appear by prolongation and elimination. However, it does not suffice to prove local existence of $\mu$ and $M$, though (\ref{movingframe}) may be rewritten as an exterior differential system in finite (but a priori unknown) dimension, for which the Cartan characters and Cartan's involutivity test may be computed (see \cite{BCG3}).

Note furthermore that a solution of (\ref{movingframe}), if it exists, does not guarantee that, among  all possible $M$, at least one is unimodular. The latter requirement, namely that $M \in \UU_{m}[\ddt]$, appears to be very restrictive, as illustrated by the non flat example of Subsection~\ref{nonflatex}.
\end{rem}

In this perspective, we give a more explicit NSC:
\begin{cor}\label{mvframinv-thm}
The $\kk[\ddt]$-ideal $\Omega$ generated by the 1-forms $\omega_{1},\ldots,\omega_{m}$ defined by (\ref{omeg-def}) is strongly closed in $\XXX_{0}$ (or, equivalently, the system $(\XX,\tau_{\XX},F)$ is flat) if, and only if, there exists $\mu \in \lin{1}{ \left( \Lambda(\XX)\right)^{m} }$, and two matrices $M\in \MM_{m,m}[\ddt]$ and $N\in \MM_{m,m}[\ddt]$ such that
\begin{equation}\label{mvframinv}
d\omega = \mu \ \omega, \qquad \dg{\mu} =  \mu^{2}, \qquad \dg{M}=- M\mu, \qquad \dg{N}= \mu N, \qquad MN=NM=I.
\end{equation}
\end{cor}
\begin{proof}
Assume that Theorem~\ref{movingframe-thm} holds true. Setting $N= M^{-1}$ and using $0=\dg{NM}=\dg{N}M+N\dg{M}$ with $\dg{M}=- M\mu$ immediately yields $\dg{N}= \mu N$ and thus (\ref{mvframinv}).
Conversely, the last equation of (\ref{mvframinv}) implies that $M\in \UU_{m}[\ddt]$, which completes the proof.
\end{proof}
Again, adapting Corollary~\ref{existmuM-cor}, this set of equations is algebraically closed.

\bigskip

From these conditions, two practical sequential procedures\footnote{We have preferred the expression ``sequential procedure'' rather than ``algorithm'' since there is no guarantee that this procedure finishes in a finite number of steps.} to test if a system is flat may be deduced. 
\paragraph{\textbf{Sequential Procedure 1.}}
We start with any meromorphic implicit system of the form (\ref{implsys}), with $P(F)$ hyper-regular.
\begin{enumerate}
\item We first compute a vector 1-form $\omega$ defined by (\ref{omeg-def}).
\item We compute the operator $\mu$, with $\ord{\mu} = \ord{\Gamma}$, such that $d\omega=\mu\omega$, given by Proposition~\ref{sol-mu}.
\item Among the possible $\mu$'s, only those satisfying $\dg{\mu} =  \mu^{2}$ are kept. If no such $\mu$ exists, we go back to step 2 and increase $\ord{\mu}$ by 1.
\item We then compute $M$ and $N$ such that $\dg{M}=- M\mu$ and $\dg{N}= \mu N$ by componentwise identification. Again, If no such meromorphic $M$ and $N$ exist, we go back to step 2, increasing the degree $\ord{\mu}$ by one.
\item Finally, only those matrices $M$ and $N$ such that $MN=NM=I$ (unimodular) are kept. If there are no such $M$'s, we go back to step 2, increasing the degree $\ord{\mu}$ by one.

If for some $\ord{\mu}$, the algorithm produces an invertible $M$, a flat output is obtained by integration of $dy=M\omega$, which is possible since $d(M\omega) =0$.
In the opposite case, the system is non flat.
\end{enumerate}

\begin{rem}
The algorithm finishes by checking if $M$ is unimodular: it suffices to apply the Smith-Jacobson algorithm that must end with the identity matrix of dimension $m$. 

Remark that, according to \cite{Ritt,Oll-aecc,KMP},  since $\ord{M^{-1}} \leq (m-1)\ord{M}$, and since $\mu=-M^{-1}\dg{M}$, we must have $\ord{\mu}\leq (m-1)\ord{M} + \ord{M}=(m)\ord{M}$, which proves that $\ord{M} \geq \frac{1}{m}\ord{\mu}$. This lower bound may be used to initialize $\ord{M}$ in step 4.

Indeed, if a solution $(\mu,M)$ exists with given $\ord{\mu}$ and $\ord{M}$ respectively, it is readily seen that a solution will exist for any $K(\mu)$ and $K(M)$ satisfying  $K(\mu)\geq \ord{\mu}$ and $K(M)\geq \ord{M}$.

Nevertheless, since we don't know upper bounds for $\ord{\mu}$ and $\ord{M}$, this algorithm remains in theory doubly infinite. It is therefore much simpler to check flatness than to prove non flatness.

Remark that, even if a system is known to be non-flat, the algorithm may provide a solution for $\mu$ and $M$ as shown by Example \ref{nonflatex}. The only obstruction in this case is that no such $M$ is unimodular.

Further details on a preliminary implementation of this sequential procedure using computer algebra may be found in \cite{AnL}.

Note again, concerning points 3 and 4, that the algorithms to find the solutions $\mu$ and $M$, if they exist, has not been addressed here.
\end{rem}
\paragraph{\textbf{Sequential Procedure 2.}}
We start again with any meromorphic implicit system of the form (\ref{implsys}), with $P(F)$ hyper-regular.
\begin{enumerate}
\item We first parameterize the matrices $M$ (of sufficiently high degree $\ord{M}$) and $N$ such that $N=M^{-1}$. More precisely, it is always possible to express the entries of $N$ as functions of the entries of $M$ to obtain such a parameterization. Note in addition that if $M$ has degree $\ord{M}$, then $N$ has at most degree $(m-1)\ord{M}$ (see again \cite{Ritt,Oll-aecc,KMP,OllBra-crm}).
\item We compute $\mu= -N\dg{M}=\dg{N}M$. Note that the latter relation always holds true since, by construction, $NM=I$, and thus $N\dg{M}+ \dg{N}M = 0$.
\item We next compute the constraints on this parameterization in order to have $d\omega=\mu\omega$.
Note that since $M=N^{-1}$, the closure condition $\dg{\mu}=\mu^{2}$ is automatically valid by Proposition~\ref{dgmu}. 
If the resulting constraints have a non void intersection, then the system is flat. Otherwise, we must go back to Step 1 and increase the degree. 
\end{enumerate}
\begin{rem}\label{discus-rem}
In a computational viewpoint, the main difficulty of the first sequential procedure relies on the need to obtain all possible solutions $\mu$, $M$ and $N$ (steps 1 to 4) in order to be sure not to miss the existence of a pair $(M,N)$ such that $M=N^{-1}$ if any.\\
For the second procedure, the main difficulties concern the parameterization of step 1 and its restriction at step 3. Let us give an example of parameterization for a $2\times 2$ unimodular matrix $M$ with $\ord{M}=1$, namely $M = M^{0} + M^{1}\ddt$. Its inverse $N$ is also, according to the bound reported in step 1, of the form $N=N^{0} + N^{1}\ddt$.
It is readily seen that a necessary and sufficient condition for $NM=MN=I$ is that
$$
\begin{array}{l}
M^{1}N^{1}= 0, \quad M^{1}N^{0} + M^{0}N^{1} + M^{1}\dot{N}^{1}= 0, \quad M^{0}N^{0} + M^{1}\dot{N}^{0}= I\\
N^{1}M^{1}= 0, \quad N^{1}M^{0} + N^{0}M^{1} + N^{1}\dot{M}^{1}= 0, \quad N^{0}M^{0} + N^{1}\dot{M}^{0}= I
\end{array}
$$
which yields\footnote{Computations done using Maple, with the help of F. Antritter} 
$$
M^{1} = m_{11}^{1}\left( \begin{array}{cc} 1&\alpha\\-\frac{1}{\beta}& -\frac{\alpha}{\beta}\end{array}\right), \quad 
N^{1} = -\gamma m_{11}^{1}\left( \begin{array}{cc} 1&\beta\\-\frac{1}{\alpha}&-\frac{\beta}{\alpha}\end{array}\right)
$$
$$
M^{0} = \left( \begin{array}{cc} m_{11}^{0}&m_{12}^{0}\\m_{21}^{0}&\frac{\alpha}{\beta}m_{11}^{0} - \frac{1}{\beta}m_{12}^{0} + \alpha m_{21}^{0}\end{array}\right), \quad 
N^{0} = \left( \begin{array}{cc} n_{11}^{0}&n_{12}^{0}\\n_{21}^{0}&n_{22}^{0}\end{array}\right)
$$
with 
$$
\begin{array}{l}
\ds n_{11}^{0} =
-\frac{\alpha}{\beta}n_{22}^{0} + \gamma \left( 2m_{11}^{0} - \frac{1}{\alpha}m_{12}^{0} + \beta m_{21}^{0} \right) + \gamma \left( \frac{\dot{\alpha}}{\alpha} + \frac{\dot{\beta}}{\beta} \right)m_{11}^{1}\\
\ds n_{12}^{0} =
-\alpha n_{22}^{0} + \gamma\beta \left( m_{11}^{0} - \frac{1}{\alpha}m_{12}^{0} \right) + \gamma\beta \frac{\dot{\alpha}}{\alpha}  m_{11}^{1}\\
\ds n_{21}^{0} = \frac{1}{\beta}n_{22}^{0} - \frac{\gamma}{\alpha} \left( m_{11}^{0} + \beta m_{21}^{0} \right) - \gamma \frac{\dot{\beta}}{\alpha\beta}  m_{11}^{1}\\
\ds m_{11}^{1} =- 
\frac
{ \left( \alpha \beta   m_{11}^{0} m_{21}^{0} n_{22}^{0}  -\beta m_{12}^{0} m_{21}^{0} n_{22}^{0}
- \beta m_{11}^{0}   + \alpha \left( m_{11}^{0}  \right) ^{2} n_{22}^{0} - m_{11}^{0} m_{12}^{0} n_{22}^{0}  \right)  \left( m_{11}^{0} +\beta m_{21}^{0} \right) 
 }
 {
 A
 }\\
 \ds \gamma=-
\frac{\alpha A}{B}
\end{array}
$$
and with
$$
\begin{array}{lcl}
A&=&\ds  - 
\dot{\beta} m_{11}^{0} +2  \dot{\alpha} \beta m_{11}^{0} m_{21}^{0} n_{22}^{0}  + 
\dot{\alpha} \left( m_{11}^{0} \right)^{2}  n_{22}^{0} 
+ \dot{\alpha} \beta^{2}  \left( m_{21}^{0} \right)^{2} n_{22}^{0}  + \beta \dot{m}_{11}^{0}  +  \beta^{2} \dot{m}_{21}^{0}
\\
\ds B &=&\ds 
 \dot{\beta} \left( m_{11}^{0} \right)^{2}m_{12}^{0}  -
 \beta  \dot{m}_{11}^{0}m_{11}^{0} m_{12}^{0}  +
  \dot{\beta} \beta m_{11}^{0} m_{12}^{0} m_{21}^{0} - 
\beta^{2}  \dot{m}_{11}^{0} m_{12}^{0} m_{21}^{0} - \beta^{2}  m_{11}^{0} m_{12}^{0} \dot{m}_{21}^{0} 
 \\
&&\ds  - 
 \beta^{3} m_{12}^{0} m_{21}^{0} \dot{m}_{21}^{0}  +
 2 \dot{\alpha} \beta^{2} \left( m_{11}^{0} \right)^{2} m_{21}^{0} 
 + \alpha \beta^{2}  \left( m_{11}^{0} \right)^{2} \dot{m}_{21}^{0} + \alpha \beta  \left( m_{11}^{0} \right)^{2} \dot{m}_{11}^{0} 
 \\
&&\ds  
+
\left( \dot{\alpha}  \beta - \alpha \dot{\beta} \right) \left( m_{11}^{0} \right)^{3}
+
\dot{\alpha} \beta^{3} m_{11}^{0} \left( m_{21}^{0} \right)^{2}  + 
\alpha \beta^{3} m_{11}^{0} m_{21}^{0} \dot{m}_{21}^{0} 
\\
&&\ds +
 \alpha \beta^{2} m_{11}^{0}\dot{m}_{11}^{0}  m_{21}^{0} - 
\alpha  \beta \dot{\beta} \left( m_{11}^{0} \right)^{2} m_{21}^{0}
 \end{array}
$$
where the 6 non zero arbitrary meromorphic functions $m_{11}^{0}$, $m_{12}^{0}$, $m_{21}^{0}$, $n_{22}^{0}$, $\alpha$ and $\beta$ constitute the required parameterization of $M$ and $N$.
The step 3 thus consists in verifying if there exists a non empty subset of the corresponding parameter space such that $d\omega=\mu\omega$, with 
$$\mu=- N\dg{M}= -(N^{0}dM^{0} + N^{1}d\dot{M}^{0}) \wedge - (N^{0}dM^{1} + N^{1}(dM^{0}+d\dot{M}^{1}))\wedge \ddt - N^{1}dM^{1}\wedge \frac{d^{2}}{dt^{2}}$$
in other words such that
$$d\omega = -(N^{0}dM^{0} + N^{1}d\dot{M}^{0})\wedge \omega - (N^{0}dM^{1} + N^{1}(dM^{0}+d\dot{M}^{1}))\wedge \dot{\omega} - N^{1}dM^{1}\wedge \ddot{\omega}.$$
In both procedures, if a unimodular matrix $M$ is found, we end up integrating the set of exterior differential equations $dy=d(M\omega)$. Formal integration algorithms are indeed far from being straightforward but their study, even restricted to our context, is beyond the scope of this paper.

To conclude this remark, let us precise that the examples presented in Section~\ref{ex-sec} below do not require using the second sequential procedure since the obtained integrating factors, at most 2 by 2 matrices, are simple enough to directly verify if they are invertible or not.
\end{rem}


\subsection{Some easy consequences}
We now show how several classical results of static feedback linearization \cite{JR,HSM}, or in the case $m=1$ \cite{CLMscl,CLM,Po-banach93,Sh,Sluis-scl93} can be recovered as consequences of Theorem~\ref{CNSthm}. 

If $\Omega$ is strongly closed, let us define $\sigma_{i}$ as the maximum degree of the entries of the $i$th column of $P(\varphi_{0})$, $i=1,\ldots, m$, which, according to Theorem~\ref{CNSthm}, locally yields:
\begin{equation}\label{flatoutmin}
x=\varphi_{0}\left(y_{1}, \ldots,y_{1}^{(\sigma_{1})},\ldots, y_{m}, \ldots,y_{m}^{(\sigma_{m})}\right).
\end{equation}
\begin{defn}
We say that a flat output $y$ is \emph{minimal} if $\sigma =\sum_{i=1}^{m}\sigma_{i}$ is minimal over all possible choices of $\widehat{U}$, $Q$ and $Z$.
\end{defn}
Obviously, a minimal $\sigma$ always exists for flat systems.
\begin{cor}\label{statlin}
A necessary and sufficient condition for system (\ref{implsys}) to be static feedback linearizable (see \cite{JR,HSM}) is that the strong closedness condition of Theorem~\ref{CNSthm} holds true and that $n=m+\sigma$ for a minimal $\sigma$.
\end{cor}
\begin{proof}
It is easily seen that (\ref{implsys}) is static feedback linearizable if, and only if, (\ref{implsys}) is \mbox{L-B equivalent} to a trivial system, the trivialization $\Phi$ being such that $\varphi_{0}$ is a local diffeomorphism, which means that $x$ is diffeomorphic to $\left(y_{1}, \ldots,y^{(\sigma_{1})},\ldots, y_{m}, \ldots,y^{(\sigma_{m})}\right)$ for a minimal $y$, i.e. $n=m+\sigma$.
\end{proof}
\begin{rem}
This Corollary generalizes the results of Jakubczyk and Respondek \cite{JR} and Hunt, Su and Meyer \cite{HSM}  in a twofold manner: first, it is not restricted to affine systems and second it applies not only in a neighborhood of an equilibrium point but of any trajectory around which the variational system is controllable. 
\end{rem}
\begin{cor}\label{singlecor}
If $m=1$, a necessary and sufficient condition for flatness is that $\Omega$ is closed in the ordinary sense, i.e. $d\omega_{1}= \tau\wedge \omega_{1}$ for some 1-form $\tau$. Furthermore, the system is flat if, and only if, it is static feedback linearizable.
\end{cor}
\begin{proof}
If $\Omega$ is generated by a single 1-form $\omega_{1}$, the unimodular matrix $M$ must be a non zero element of $\kk$ and the strong closedness of $\Omega$ reduces to ordinary closedness. Thus the system is flat if, and only if, $\Omega$ is closed. Assuming closedness,  by Frobenius' Theorem, we deduce that there exists a scalar function $\psi_{0}$ such that $d\psi_{0}=M\omega_{1}$ and that (\ref{flatoutmin}) holds true with
$\sigma=n-1$. Using Corollary~\ref{statlin}, the system is static feedback linearizable. The converse is trivial since every static feedback linearizable system is flat.
\end{proof}


\section{Examples}\label{ex-sec}
\subsection{Non holonomic car}
Consider the 3 dimensional system in the $x-y$ plane, representing a vehicle of length $l$, whose orientation is given by the angle $\theta$, the coordinates $(x,y)$ standing for the position of the middle of the rear axle, and controlled by the velocity modulus $u$ and the angular position of the front wheels $\varphi$.
\begin{equation}\label{nonholcar}
\begin{array}{l}
\dot{x}=u\cos\theta\\
\dot{y}=u\sin\theta\\
\dot{\theta}=\frac{u}{l}\tan\varphi
\end{array}
\end{equation}
Since $n=3$ and $m=2$, $n-m=1$ and (\ref{nonholcar}) is equivalent to the single implicit equation obtained by eliminating the inputs $u$ and $\varphi$:
\begin{equation}\label{implnonholcar}
F(x,y,\theta,\dot{x},\dot{y},\dot{\theta})=\dot{x}\sin\theta - \dot{y}\cos\theta =0
\end{equation}
We immediately have:
\begin{equation}\label{PFcar}
\begin{array}{l}
\ds P(F)= \left( \frac{\partial F}{\partial x}+\frac{\partial F}{\partial \dot{x}}\ddt, \quad \frac{\partial F}{\partial y}+\frac{\partial F}{\partial \dot{y}}\ddt, \quad
\frac{\partial F}{\partial \theta}+\frac{\partial F}{\partial \dot{\theta}}\ddt \right)\\
\hspace{2cm}\ds = \left( \sin\theta \ddt, \quad -\cos\theta \ddt, \quad \dot{x}\cos\theta + \dot{y}\sin\theta \right).
\end{array}
\end{equation}
Setting $E=\dot{x}\cos\theta + \dot{y}\sin\theta$, we apply the Smith decomposition algorithm of Appendix~\ref{Smith-app}: moving the last column (of degree zero) to the first place by a permutation with the two others, we get $P(F)U_0=\left( E, \quad -\cos\theta\ddt, \quad \sin\theta\ddt \right)$ with $U_0 \triangleq \left(\begin{array}{ccc} 0&0&1\\0&1&0\\1&0&0\end{array}\right)$, and then, again right-multiplying the result by $U_1 \triangleq  \left(\begin{array}{ccc}\frac{1}{E}&\frac{\cos\theta}{E}\ddt&-\frac{\sin\theta}{E}\ddt\\0&1&0\\0&0&1
\end{array}\right)$ yields
$$P(F)U=\left( 1,\quad 0,\quad 0 \right), \quad \mbox{\textrm{with}~} U \triangleq U_0 U_1 =\left( \begin{array}{ccc}0&0&1\\0&1&0\\\frac{1}{E}&\frac{\cos\theta}{E}\ddt&-\frac{\sin\theta}{E}\ddt \end{array}\right).$$
Thus $P(F)$ is hyper-regular and
$$\widehat{U}=U\left(\begin{array}{c}0_{1,2}\\I_{2}\end{array}\right) =\left( \begin{array}{cc}0&1\\1&0\\\frac{\cos\theta}{E}\ddt&-\frac{\sin\theta}{E}\ddt \end{array}\right)$$
with $I_{2}$ the identity matrix of $\RR^{2}$.
Computing $Q_0$ by (\ref{Q0def-Id}) yields
$$Q_0 =\left( \begin{array}{ccc}0&1&0\\1&0&0\\\frac{\sin\theta}{E}\ddt&-\frac{\cos\theta}{E}\ddt&1  \end{array}\right), \quad \widetilde{Q}_0 = \left( \begin{array}{ccc}0&1&0\\1&0&0  \end{array}\right).$$

Multiplying $Q_0$ by the vector $\left( \begin{array}{c}dx\\dy\\d\theta\end{array}\right)$,  the last line  reads
$\frac{1}{E}\left( \sin\theta d\dot{x}-\cos\theta d\dot{y}+(\dot{x}\cos\theta +\dot{y}\sin\theta)d\theta \right)=\frac{1}{E} d(\dot{x}\sin\theta - \dot{y}\cos\theta)$
and, by (\ref{implnonholcar}), identically vanishes on $\XX_{0}$.

Setting $\widetilde{Q}_0 \left( \begin{array}{c}dx\\dy\\d\theta\end{array}\right) =\left( \begin{array}{c}\omega_{1}\\\omega_{2}\end{array}\right)$, i;e. $\omega_{1}=dy$ and $\omega_{2}=dx$, the ideal $\Omega$ generated by $(\omega_{1}, \omega_{2})$ is trivially strongly closed with $M=I_{2}$, which finally gives the flat output
$y_{1} = y$ and $y_{2} = x$. 
We have thus recovered the flat output originally obtained in \cite{rouchon-et-al-cdc93,rouchon-et-al-ecc93}, up to a permutation of the components of $y$.
\subsection{Non holonomic car (continued)}

Other decompositions of  $P(F)$, given by (\ref{PFcar}), may indeed be obtained, but they are all equivalent in the sense that one decomposition may be deduced from another one by multiplication by a unimodular matrix. However, the resulting vector 1-form $\omega$, contrarily to what happens in the previous example, may not be integrable. Our aim is here to show how the generalized moving frame structure equations (\ref{movingframe}) may be used to obtain an integrable $M\omega$. Such an example is provided by restarting the right-Smith decomposition of  $P(F)$ by right-multiplying it by $\left(\begin{array}{ccc}\cos\theta&0&0\\\sin\theta&1&0\\0&0&1\end{array}\right)$ and using the formula $\sin\theta\ddt (\cos\theta) -\cos\theta\ddt(\sin\theta)=-\dot{\theta}$, we obtain
$$
U=\left(\begin{array}{ccc}
\cos\theta&-\frac{1}{\dot{\theta}}\cos^{2}\theta\ddt&\frac{1}{\dot{\theta}}\left( \dot{x}\cos\theta+\dot{y}\sin\theta\right)\cos\theta\\
\sin\theta&1-\frac{1}{\dot{\theta}}\sin\theta\cos\theta\ddt&\frac{1}{\dot{\theta}}\left(
\dot{x}\cos\theta+\dot{y}\sin\theta\right)\sin\theta\\
0&0&1
\end{array}\right).
$$
Using (\ref{implnonholcar}), we get $\dot{x}\cos\theta+\dot{y}\sin\theta = \frac{\dot{x}}{\cos\theta}$, and 
$\frac{1}{\dot{\theta}}\left( \dot{x}\cos\theta+\dot{y}\sin\theta\right)\cos\theta = \frac{\dot{x}}{\dot{\theta}}$, 
$\frac{1}{\dot{\theta}}\left( \dot{x}\cos\theta+\dot{y}\sin\theta\right)\sin\theta = \frac{\dot{y}}{\dot{\theta}}$. Thus
$$
U=\left(\begin{array}{ccc}
\cos\theta&-\frac{1}{\dot{\theta}}\cos^{2}\theta\ddt&\ds \frac{\dot{x}}{\dot{\theta}}\\
\sin\theta&1-\frac{1}{\dot{\theta}}\sin\theta\cos\theta\ddt&\ds \frac{\dot{y}}{\dot{\theta}}\\
0&0&1
\end{array}\right).
$$
and
$$
\widehat{U}=\left(\begin{array}{cc}
-\frac{1}{\dot{\theta}}\cos^{2}\theta\ddt&\ds \frac{\dot{x}}{\dot{\theta}}\\
1-\frac{1}{\dot{\theta}}\sin\theta\cos\theta\ddt&\ds \frac{\dot{y}}{\dot{\theta}}\\
0&1
\end{array}\right)
$$
and then, applying (\ref{Q0def-Id}):
$$Q_0=\left(\begin{array}{ccc}
-\tan\theta&1&0\\
0&0&1\\
-\frac{1}{\dot{\theta}}\sin\theta\cos\theta\ddt&\frac{1}{\dot{\theta}}\cos^{2}\theta\ddt&
\ds -\frac{\dot{x}}{\dot{\theta}}
\end{array}\right).
$$

The vector 1-form $\omega= (\omega_{1},\omega_{2})^{T}$ is again obtained by multiplying the two first rows of $Q_0$ by  $(dx,dy,d\theta)^{T}$:
$\omega_{1}=-\tan\theta dx+dy$ and $\omega_{2}=d\theta$. 

We have $d\omega_{1}=-\frac{1}{\cos^{2}\theta}d\theta\wedge dx=\frac{1}{\dot{\theta}}\omega_{2}\wedge \dot{\omega}_{1}$ and thus
$\left(\begin{array}{c}d\omega_{1}\\d\omega_{2}\end{array}\right) =
\left(\begin{array}{c}\frac{1}{\dot{\theta}}\omega_{2}\wedge \dot{\omega}_{1}\\0\end{array}\right) =
\left(\begin{array}{cc}0&-\frac{1}{\dot{\theta}} \dot{\omega}_{1}\wedge \\0&0\end{array}\right) 
\left(\begin{array}{c}\omega_{1}\\\omega_{2}\end{array}\right)
\triangleq \left(\begin{array}{cc}0&\Gamma_{1,1,0}^{1,2} \dot{\omega}_{1}\wedge \\0&0\end{array}\right) \omega$, with $\ord{\Gamma}=0$.

Hence, according to step 2 of the Sequential Procedure 1, a possible choice is $\mu= \left( \begin{array}{cc}0&\left( \Gamma_{1,1,0}^{1,2} \dot{\omega}_{1}+ \nu_{1,0,0}^{2,2}\omega_{2}\right)\wedge \\0&0\end{array}\right)$ where $\nu_{1,0,0}^{2,2}$ is an arbitrary meromorphic function. Since we have $\mu^{2}=\left( \begin{array}{cc}0&0\\0&0\end{array}\right)$, we must have $d\left( \Gamma_{1,1,0}^{1,2} \dot{\omega}_{1}+ \nu_{1,0,0}^{2,2}\omega_{2}\right)=0$. Since $\left(dx,dy,d\theta\right)^{T}=\widehat{U}\omega$ by Corollary~\ref{flatvar-cor}, we have
 $dx= -\frac{1}{\dot{\theta}}\cos^{2}\theta \dot{\omega}_{1} + \frac{\dot{x}}{\dot{\theta}}\omega_{2}$ and thus $\Gamma_{1,1,0}^{1,2}\dot{\omega}_{1}= -\frac{1}{\dot{\theta}} \dot{\omega}_{1}= \frac{1}{\cos^{2}\theta}dx - \frac{\dot{x}}{\dot{\theta}\cos^{2}\theta}\omega_{2}= \frac{1}{\cos^{2}\theta}dx - \frac{\dot{x}}{\dot{\theta}\cos^{2}\theta}d\theta$. Thus, choosing $\nu_{1,0,0}^{2,2}= \frac{\dot{x}}{\dot{\theta}\cos^{2}\theta} + 2\frac{x\sin\theta}{\cos^{3}\theta}$, we get
 $\Gamma_{1,1,0}^{1,2} \dot{\omega}_{1}+ \nu_{1,0,0}^{2,2}\omega_{2}= \frac{1}{\cos^{2}\theta}dx + 2\frac{x\sin\theta}{\cos^{3}\theta} d\theta = d\left(\frac{x}{\cos^{2}\theta}\right)$ and therefore, as desired, $\dg{\mu}=0=\mu^{2}$ with 
 $$\mu= \left( \begin{array}{cc}0&d\left(\frac{x}{\cos^{2}\theta}\right)\wedge \\0&0\end{array}\right)$$
 
 Again, by componentwise identification of $\dg{M}=-M\mu$, choosing $M$ of the form $M= \left(\begin{array}{cc} 1&\ast\\0&1\end{array}\right)$, with $\ord{M}=0$, $M$ being thus unimodular by construction, one immediately finds $M=\left(\begin{array}{cc} 1&-\frac{x}{\cos^{2}\theta}\\0&1\end{array}\right)$.
and then
$$M\omega= \left(\begin{array}{c}-\tan\theta dx - \frac{x}{\cos^{2}\theta}d\theta +dy\\d\theta\end{array}\right),\quad d(M\omega)=0.$$
Thus, setting $\left(\begin{array}{c}dy_{1}\\dy_{2}\end{array}\right)=M\omega$, we obtain
$$y_{1}=y-x\tan\theta,\quad y_{2}=\theta$$
which is another possible flat output: it is easily checked that the inverse L-B isomorphism is given by $x=-\frac{\dot{y}_{1}}{\dot{y}_{2}}\cos^{2}y_{2}$, $y=y_{1}-\frac{\dot{y}_{1}}{\dot{y}_{2}}\sin y_{2}\cos y_{2}$, $\theta=y_{2}$.
\subsection{The pendulum}
We consider a pendulum in the vertical plane studied in \cite{FLMR-ieee},  of length $l$ and inertia $J$, whose mass $m$ is concentrated at its end point $C$. It is  controlled by the two components of the force $F$ applied to the opposite end point $A$ of the pendulum. Introducing an inertial frame $(0,x,z)$, it is modelled by 
\begin{equation}\label{pendul-eq}
\begin{array}{l}
\ds \ddot{x}=u_{1}\\
\ds \ddot{z}=u_{2}\\
\ds a\ddot{\theta}= -u_{1}\cos\theta + (u_{2}+1)\sin\theta
\end{array}
\end{equation}
with $x=\frac{x_{C}}{g}$, $z=\frac{z_{C}}{g}$, $(x_{C},z_{C})$ being the coordinates of $C$, $\theta$ being the angle between the pendulum and the vertical axis, $u_{1}=\frac{F_{x}}{mg}$, $u_{2}=\frac{F_{z}}{mg}-1$, $(F_{x},F_{z})$ being the components of the force $F$, and $a=\frac{J}{mgl}$.

An implicit model is given by:
\begin{equation}\label{pendul-impl}
F(x,\dot{x},\ddot{x},z,\dot{z},\ddot{z},\theta,\dot{\theta}, \ddot{\theta})=a\ddot{\theta} + \ddot{x}\cos\theta - ( \ddot{z}+1)\sin\theta = 0.
\end{equation}
Though this is a second order system, it can be easily transformed into a first order one by setting $\dot{x}=v_{x}$, $\dot{z}=v_{z}$ and $\dot{\theta}=v_{\theta}$. We thus obtain the 4 dimensional implicit system
 \begin{equation}\label{pendul-1impl}
 \begin{array}{l}
 \ds \dot{x}-v_{x}=0\\
 \ds \dot{z}-v_{z}=0\\
 \ds \dot{\theta}-v_{\theta}=0\\
 \ds a\dot{v}_{\theta} + \dot{v}_{x}\cos\theta - ( \dot{v}_{z}+1)\sin\theta = 0.
 \end{array}
 \end{equation}
However, it can be easily verified that all the results of this paper can be extended word for word to higher order systems. Thus, because of its smaller dimension, we prefer using (\ref{pendul-impl}) instead of (\ref{pendul-1impl}).

The variational system corresponding to (\ref{pendul-impl}) is given by
\begin{equation}\label{cotpend}
\left( \cos\theta\frac{d^{2}}{dt^{2}}, -\sin\theta\frac{d^{2}}{dt^{2}}, a\frac{d^{2}}{dt^{2}}-b \right)\left(\begin{array}{c}dx\\dz\\d\theta\end{array}\right)=0
\end{equation}
where $b=\ddot{x}\sin\theta +(\ddot{z}+1)\cos\theta$.

Using the identity $\left(\cos\theta \frac{d^{2}}{dt^{2}}\right)\left( -a\cos\theta \right) -\left(\sin\theta\frac{d^{2}}{dt^{2}}\right) \left( a\sin\theta \right) + \left( a\frac{d^{2}}{dt^{2}} -b\right)(1) = a\dot{\theta}^{2}-b$ and setting $E=a\dot{\theta}^{2}-b$, we have
$$\left( \cos\theta\frac{d^{2}}{dt^{2}},  -\sin\theta\frac{d^{2}}{dt^{2}}, a\frac{d^{2}}{dt^{2}}-b \right)
\left( \begin{array}{ccc}-a\cos\theta&0&1\\a\sin\theta&1&0\\1&0&0\end{array}\right)
\left( \begin{array}{ccc}\frac{1}{E}&\frac{\sin\theta}{E}\frac{d^{2}}{dt^{2}}&-\frac{\cos\theta}{E}\frac{d^{2}}{dt^{2}}\\0&1&0\\0&0&1\end{array}\right) =
\left( 1,0,0\right)$$
which proves the hyper-regularity of $P(F)$ and we get
$$ U=\left( \begin{array}{ccc}
-\frac{a\cos\theta}{E}&-\frac{a\sin\theta\cos\theta}{E}\frac{d^{2}}{dt^{2}}&\frac{a\cos^{2}\theta}{E}\frac{d^{2}}{dt^{2}}+1\\
\frac{a\sin\theta}{E}&
\frac{a\sin^{2}\theta}{E}\frac{d^{2}}{dt^{2}}+1&-\frac{a\sin\theta\cos\theta}{E}\frac{d^{2}}{dt^{2}}\\
\frac{1}{E}&
\frac{\sin\theta}{E}\frac{d^{2}}{dt^{2}}&-\frac{\cos\theta}{E}\frac{d^{2}}{dt^{2}}
\end{array}\right),\quad 
\widehat{U}=\left( \begin{array}{cc}
-\frac{a\sin\theta\cos\theta}{E}\frac{d^{2}}{dt^{2}}&\frac{a\cos^{2}\theta}{E}\frac{d^{2}}{dt^{2}}+1\\
\frac{a\sin^{2}\theta}{E}\frac{d^{2}}{dt^{2}}+1&-\frac{a\sin\theta\cos\theta}{E}\frac{d^{2}}{dt^{2}}\\
\frac{\sin\theta}{E}\frac{d^{2}}{dt^{2}}&-\frac{\cos\theta}{E}\frac{d^{2}}{dt^{2}}
\end{array}\right).
$$
Now left decomposing $\widehat{U}$, we get
$$
\left( \begin{array}{ccc}
1&0&0\\0&1&0\\-\frac{\sin\theta}{E}\frac{d^{2}}{dt^{2}}&\frac{\cos\theta}{E}\frac{d^{2}}{dt^{2}}&1
\end{array}\right)
\left( \begin{array}{ccc}
0&1&-a\sin\theta\\1&0&a\cos\theta\\0&0&1
\end{array}\right)
\widehat{U}=
\left(\begin{array}{cc}1&0\\0&1\\0&0\end{array}\right).
$$
Hence, 
$$Q_0=\left(\begin{array}{ccc}
0&1&-a\sin\theta\\1&0&a\cos\theta\\\frac{\cos\theta}{E}\frac{d^{2}}{dt^{2}}&-\frac{\sin\theta}{E}\frac{d^{2}}{dt^{2}}&\frac{a}{E}\frac{d^{2}}{dt^{2}}-\frac{b}{E}
\end{array}\right)$$
and
$$Q_0\left(\begin{array}{c}dx\\dz\\d\theta\end{array}\right)=
\left(\begin{array}{c}
dz-a\sin\theta d\theta\\dx+a\cos\theta d\theta\\\frac{1}{E}\left( \cos\theta d\ddot{x}-\sin\theta d\ddot{z} +a d\ddot{\theta}-b d\theta\right)
\end{array}\right)
=
\left(\begin{array}{c}
dz-a\sin\theta d\theta\\dx+a\cos\theta d\theta\\0
\end{array}\right)
=
\left(\begin{array}{c}
d(z+a\cos\theta)\\d(x+a\sin\theta)\\0
\end{array}\right)
$$
Thus, the strong closedness condition holds true: setting $M=I_{2}$, we obtain
$$d(z+a\cos\theta)=dy_{1},\quad d(x+a\sin\theta)=dy_{2}$$
or
$$y_{1}=z+a\cos\theta, \quad y_{2}=x+a\sin\theta$$
which represents, up to a permutation of $y_{1}$ and $y_{2}$, the coordinates of the Huygens oscillation center, already found in \cite{pM1,FLMR-ieee}.
\subsection{A non flat example}\label{nonflatex}

We consider the (single-input) system with the coordinates $x_1$ and $x_2$ in implicit form\protect{\footnote{the author is indebted to Dr. Felix Antritter for his help in the computations of this example in \textit{Maple}.}}
\begin{equation}\label{impoliv}
\dot{x}_2-\frac{1}{2}(\dot{x}_1)^2=0\,.
\end{equation}
This implicit system corresponds to the explicit one
$$
\dot{x}_1=u,\quad \dot{x}_2=\frac{1}{2}u^2
$$
which is notably non-static feedback linearizable, since according to \cite{JR,HSM} we set
$f=u\frac{\partial}{\partial x_1}+\frac{1}{2}u^2\frac{\partial}{\partial x_2}$ and $g=\frac{\partial}{\partial u}$ and compute the following distributions: $G_0=\span{g}=\ol{G}_0$ (where $\ol{G}_0$ is the involutive closure of $G_0$), $G_1=G_0+[f,G_0]=\span{\frac{\partial}{\partial u},(-\frac{\partial}{\partial x_1}-u\frac{\partial}{\partial x_2})}\neq \ol{G}_1=T\RR^3$. We conclude that the system is not flat according to the equivalence proven in \cite{CLMscl,CLM} that for single-input systems static and dynamic feedback linearization are equivalent.

It can also be verified that this system is not flat since it doesn't satisfy the ruled manifold criterion \cite{R1,Sluis-scl93}:
for all $(x_1,x_2,\dot{x}_1,\dot{x}_2)$ satisfying (\ref{impoliv}), there must exist a non zero vector $(g_1,g_2)$ such that
\begin{equation}
(\dot{x}_2+\lambda g_2)-\frac{1}{2}(\dot{x}_1+\lambda g_1)^2\equiv0
\end{equation}
holds for all $\lambda$ in an open interval of $\RR$ containing the origin. Developing this expression, we get the following second degree polynomial with respect to $\lambda$:
$$(\dot{x}_2-\frac{1}{2}(\dot{x}_1)^2) +(g_2-\dot{x}_1g_1)\lambda -\frac{1}{2}(g_1)^2\lambda^2 \equiv0$$
Since the 0th order term is equal to 0 according to (\ref{impoliv}, the coefficients of $\lambda^2$ and $\lambda$ have to vanish identically, i.e. $g_1=0$, and thus $g_2=0$, which proves that the ruled manifold criterion is not satisfied.

We finally verify non flatness using the Generalized Moving Frame Structure Equations.

The variational system is
\begin{equation}\label{tanoliv}
P(F)dx=\left( \begin{array}{cc} -\dot{x}_{1}\ddt &\ddt \end{array}\right) \left( \begin{array}{c} dx_{1}\\dx_{2}\end{array}\right) = 0.
\end{equation}
We compute $\widehat{U}$ given by (\ref{Uhatdef}) with $U\in \rsm{P(F)}$:
\begin{equation}
U=
\left( \begin{array}{cc}
\frac{1}{\ddot{x}_{1}}&- \frac{1}{\ddot{x}_{1}}\ddt\\
\frac{\dot{x}_{1}}{\ddot{x}_{1}}&1- \frac{\dot{x}_{1}}{\ddot{x}_{1}}\ddt \end{array}\right), \,\,\quad
\widehat{U}= \left( \begin{array}{c}
- \frac{1}{\ddot{x}_{1}}\ddt\\
1- \frac{\dot{x}_{1}}{\ddot{x}_{1}}\ddt 
\end{array}\right)
\end{equation}
Then
\begin{equation}
Q_0=\left( \begin{array}{cc} 
-\dot{x}_{1}&1\\
- \frac{\dot{x}_{1}}{\ddot{x}_{1}}\ddt & \frac{1}{\ddot{x}_{1}}\ddt \end{array}\right),\;\quad
\widetilde{Q}_0=\left( 
\begin{array}{cc} 
-\dot{x}_{1}&1
\end{array}
\right).
\end{equation}
Thus, a flat output $\omega=\widetilde{Q}(\,dx_1\;dx_2\,)^T$ of the variational system is obtained as
\begin{equation}\label{omegaoliv}
\omega=-\dot{x}_1dx_1+dx_2
\end{equation}
with, according to Corollary~\ref{flatvar-cor},
\begin{equation}\label{dx1omega}
dx_1=-\frac{1}{\ddot{x}_1}\dot{\omega},\quad dx_2= \omega -\frac{\dot{x}_{1}}{\ddot{x}_{1}} \dot{\omega}\,.
\end{equation}
The exterior derivative of $\omega$ is non zero
\begin{equation}\label{domegaoliv}
d\omega= - d\dot{x}_{1}\wedge dx_{1}= -\frac{1}{\ddot{x}_{1}^{2}}\ddot{\omega}\wedge \dot{\omega}.
\end{equation}
We compute $\mu$ satisfying the first two equations of the Generalized Moving Frame Structure Equations. 

From Proposition~\ref{sol-mu}, the general solution of $d\omega=\mu\omega$ is given by
\begin{equation}\label{exnfmu}
\mu=-\frac{1}{(\ddot{x}_1)^2}\ddot{\omega}\wedge\ddt+\sum_{i,j\geq0}\mu_{i,j}\omega^{(i)}\wedge\frac{d^j}{dt^j},\quad \mu_{i,j}-\mu_{j,i}=0\quad \forall i>j\,,
\end{equation}
where the $\mu_{i,j}$ are meromorphic functions of $\ol{x}$.

The first term in the right hand side of (\ref{exnfmu}) being of the first order, the minimum degree of the remaining polynomial must be one. We show that $\mu$ satisfying $\dg{\mu} =  \mu^{2}$ can be obtained in the following particular form of degree one:
$$\mu=-\frac{1}{(\ddot{x}_1)^2}\ddot{\omega}\wedge\ddt+\mu_{00}\omega\wedge+\mu_{11}\dot{\omega}\wedge\ddt.$$ 
We compute
\begin{eqnarray*}
-\dg{\mu}
&=&
\left[-d\mu_{00} \wedge \omega
 - \mu_{00} d\omega\right]\wedge
\\&&
+
\left[- \left(2 {\ddot{x}_{1}}^{-3}\right) d\ddot{x}_{1} \wedge \ddot{\omega}
 + \left({\ddot{x}_{1}}^{-2}\right) d\ddot{\omega}
 - d\mu_{11} \wedge \dot{\omega}
 - \mu_{11} d\dot{\omega}\right]\wedge \frac{d}{dt}
\end{eqnarray*}
After substitution of $d\ddot{x}_1$ using
\begin{equation}\label{dd2tx1}
d\ddot{x}_1=
\left(-2 {\frac {{x_{1}^{(3)}}^2 }{{\ddot{x}_{1}}^3 }}+{\frac {x_{1}^{(4)}}{{\ddot{x}_{1}}^2 }}\right)\cdot \dot{\omega}
 + \left(2 {\frac {x_{1}^{(3)}}{{\ddot{x}_{1}}^2 }}\right)\cdot \ddot{\omega}
 + \left(-{\ddot{x}_{1}}^{-1}\right)\cdot \omega^{(3)}
\end{equation}
and time derivatives of $d\omega$, we get
\begin{eqnarray*}
-\dg{\mu}&=&
 \left[-d\mu_{00} \wedge \omega
 + \left({\frac {\mu_{00}}{{\ddot{x}_{1}}^2 }}\right) \ddot{\omega} \wedge \dot{\omega}\right]\wedge
\\&&
+\left[
\left(2 {\frac {5 {x_{1}^{(3)}}^2 -2 x_{1}^{(4)}\ddot{x}_{1}+\mu_{11}x_{1}^{(3)}{\ddot{x}_{1}}^3 }{{\ddot{x}_{1}}^6 }}\right) \dot{\omega} \wedge \ddot{\omega}
 + \left({\ddot{x}_{1}}^{-4}\right) \omega^{(3)} \wedge \ddot{\omega}
\right]
\wedge \frac{d}{dt}
\\&&
 +
\left[
 \left({\frac {4 x_{1}^{(3)}+\mu_{11}{\ddot{x}_{1}}^3 }{{\ddot{x}_{1}}^5 }}\right)\omega^{(3)} \wedge \dot{\omega}
 + \left(-{\ddot{x}_{1}}^{-4}\right) \omega^{(4)} \wedge \dot{\omega}
- d\mu_{11} \wedge \dot{\omega}
\right]\wedge\frac{d}{dt}
\end{eqnarray*}
The zeroth order term of $\mu^2$ is
\begin{eqnarray*}
\left[ - \left({\frac {\dot{\mu}_{00}}{{\ddot{x}_{1}}^2 }}\right) \ddot{\omega} \wedge \omega
 - \left({\frac {\mu_{00}}{{\ddot{x}_{1}}^2 }}\right) \ddot{\omega} \wedge \dot{\omega}
 + \left(  \mu_{11}\dot{\mu}_{00}\right)  \dot{\omega} \wedge \omega
\right]\wedge 
\end{eqnarray*}
and the first order term is
\begin{eqnarray*}
&&
\left[
\left({\ddot{x}_{1}}^{-4}\right) \ddot{\omega} \wedge \omega^{(3)}
 - \left({\frac {\dot{\mu}_{11}\ddot{x}_{1}+2 \mu_{11}x_{1}^{(3)}+{\mu_{11}}^2 {\ddot{x}_{1}}^3 }{{\ddot{x}_{1}}^3 }}\right) \ddot{\omega} \wedge \dot{\omega}
 - \left({\frac {\mu_{11}}{{\ddot{x}_{1}}^2 }}\right) \dot{\omega} \wedge \omega^{(3)}
\right]\wedge \frac{d}{dt}
\end{eqnarray*}
Note that the second order term of $\mu^2$ identically vanishes. 

Combining the previous formulas, the zeroth order term of $\mu^2-\dg{\mu}$ is
\begin{eqnarray*}
\left[ -
\left({\frac {\dot{\mu}_{00}}{{\ddot{x}_{1}}^2 }}\right) \ddot{\omega} \wedge \omega
 + \left(\mu_{11}\dot{\mu}_{00}\right) \dot{\omega} \wedge \omega
 - d\mu_{00} \wedge \omega
\right]\wedge
\end{eqnarray*}
and the first order term is
\begin{eqnarray*}
&&
-\left({\frac {{\ddot{x}_{1}}^4 \dot{\mu}_{11}+4 \mu_{11}x_{1}^{(3)}{\ddot{x}_{1}}^3 +{\ddot{x}_{1}}^6 {\mu_{11}}^2 +10 {x_{1}^{(3)}}^2 -4 x_{1}^{(4)}\ddot{x}_{1}}{{\ddot{x}_{1}}^6 }}\right)\ddot{\omega} \wedge \dot{\omega}\wedge \frac{d}{dt}
\\&&
-
\left[
2\left( {\frac {\mu_{11}{\ddot{x}_{1}}^3 +2 x_{1}^{(3)}}{{\ddot{x}_{1}}^5 }}\right) \dot{\omega} \wedge \omega^{(3)}
 + \left({\ddot{x}_{1}}^{-4}\right) \omega^{(4)} \wedge \dot{\omega}
 + d\mu_{11} \wedge \dot{\omega}
\right]\wedge \frac{d}{dt}
\end{eqnarray*}
Assuming that $\mu_{11}$ depends on $(x_1,x_2,\dot{x}_1,\ddot{x}_1,{x}^{(3)}_1)$, the above equation, with $d\mu_{11}=\sum_{i=0}^3\frac{\partial \mu_{11}}{\partial x_1^{(i)}}dx_1^{(i)}+\frac{\partial \mu_{11}}{\partial x_2}dx_2$ and after substituting the successive derivatives of $dx_{1}$ in function of $\omega$ and derivatives, reads
\begin{eqnarray*}
&&
-\left({\frac {E}{\left( \ddot{x}_{1}\right)^6 }}\right)\ddot{\omega} \wedge \dot{\omega}\wedge \frac{d}{dt}
\\&&
 - \left({\frac {2 \mu_{11}\left( \ddot{x}_{1}\right)^3 +4 x_{1}^{(3)}+ \frac{\partial \mu_{11}}{\partial \ddot{x}_1} \left(\ddot{x}_{1}\right)^4 -3  \frac{\partial \mu_{11}}{\partial x_1^{(3)}} \left(\ddot{x}_{1}\right)^3 x_{1}^{(3)}}{\left(\ddot{x}_{1}\right)^5 }}\right) \dot{\omega} \wedge \omega^{(3)}\wedge \frac{d}{dt}
\\&&
 -
\left[
\left({\frac {1- \frac{\partial \mu_{11}}{\partial x_1^{(3)}} \left(\ddot{x}_{1}\right)^3 }{\left(\ddot{x}_{1}\right)^4 }}\right) \omega^{(4)} \wedge \dot{\omega}
 - \left( \frac{\partial \mu_{11}}{\partial {x}_2} \right)\omega \wedge \dot{\omega}
\right]\wedge \frac{d}{dt}
\end{eqnarray*}
with the expression
$$
\begin{array}{l}
E\triangleq
\left(\ddot{x}_{1}\right)^4 \dot{\mu}_{11}+4 \mu_{11}x_{1}^{(3)}\left(\ddot{x}_{1}\right)^3 +\left(\ddot{x}_{1}\right)^6 \left(\mu_{11}\right)^2 +10 \left(x_{1}^{(3)}\right)^2 -4 x_{1}^{(4)}\ddot{x}_{1}- \frac{\partial \mu_{11}}{\partial \dot{x}_1} \left(\ddot{x}_{1}\right)^5\\
\hspace{1.5cm} +2  \frac{\partial \mu_{11}}{\partial \ddot{x}_1} x_{1}^{(3)}\left(\ddot{x}_{1}\right)^4 -6  \frac{\partial \mu_{11}}{\partial x_1^{(3)}} \left(\ddot{x}_{1}\right)^3 \left(x_{1}^{(3)}\right)^2 +3  \frac{\partial \mu_{11}}{\partial x_1^{(3)}} \left(\ddot{x}_{1}\right)^4 x_{1}^{(4)}.
\end{array}
$$

After substitution of $\dot{\mu}_{11}$ by $\sum_{i=0}^3\frac{\partial \mu_{11}}{\partial x_1^{(1)}}x_1^{(i+1)}+\frac{\partial \mu_{11}}{\partial x_2}\dot{x}_2$, and using the system equation, we get the PDEs
\begin{eqnarray*}
-{\frac {\partial \mu_{11}}{\partial x_{1}}} \left(\ddot{x}_1\right)^4\dot{x}_1-\frac{1}{2} {\frac {\partial \mu_{11}}{\partial x_{2}}}\left(\ddot{x}_1\right)^4 \left(\dot{x}_1\right)^2 -3 {\frac {\partial \mu_{11}}{\partial (\ddot{x}_1)}}{x_1^{(3)}} \left(\ddot{x}_1\right)^4
 -4 {\frac {\partial \mu_{11}}{\partial {x_1^{(3)}}}} \left(\ddot{x}_1\right)^4{x_1^{(4)}}
&&\\
-4 \mu_{11} {x_1^{(3)}} \left(\ddot{x}_1\right)^3 -\left(\mu_{11}\right)^2 \left(\ddot{x}_1\right)^6
 -10 \left({x_1^{(3)}}\right)^2 +4 {x_1^{(4)}}\ddot{x}_1+6 {\frac {\partial \mu_{11}}{\partial {x_1^{(3)}}}} \left(\ddot{x}_1\right)^3 \left({x_1^{(3)}}\right)^2
&=&
0\\
-2 \mu_{11} \left(\ddot{x}_1\right)^3 -4{x_1^{(3)}}-{\frac {\partial \mu_{11}}{\partial \ddot{x}_1}} \left(\ddot{x}_1\right)^4 +3 {\frac {\partial \mu_{11}}{\partial {x_1^{(3)}}}} \left(\ddot{x}_1\right)^3 {x_1^{(3)}}
&=&
0\\
-1+{\frac {\partial \mu_{11}}{\partial {x_1^{(3)}}}} \left(\ddot{x}_1\right)^3 
&=&
0\\
-{\frac {\partial \mu_{11}}{\partial x_{2}}}
&=&
0
\end{eqnarray*}
One can verify that this set of PDE's admits the following solution\footnote{thanks to\textit{Maple}}:
\begin{equation}
\mu_{11}
=
{\frac {x_{1}^{(3)}}{\left(\ddot{x}_{1}\right)^3 }}+{\frac {\dot{x}_{1}}{(x_{1}+F_1(\dot{x}_{1}) \dot{x}_{1}) \left(\ddot{x}_{1}\right)^2 }}
\end{equation}
where $F_1(\dot{x}_{1})$ is an arbitrary meromorphic function of $\dot{x}_{1}$.
For simplicity's sake, we choose  $F_1= 0$:
\begin{equation}\label{mu11sol}
\mu_{11}
=
{\frac {x_{1}^{(3)}}{{\ddot{x}_{1}}^3 }}+{\frac {\dot{x}_{1}}{x_{1} {\ddot{x}_{1}}^2 }}
\end{equation}
thus implying the following zeroth order term in $\mu^2-\dg{\mu}$:
$$
\left[
\left(-{\frac {\dot{\mu}_{00}}{{\ddot{x}_{1}}^2 }}\right)\ddot{\omega} \wedge \omega
 + \left({\frac {(x_{1}^{(3)}x_{1}+\dot{x}_{1}\ddot{x}_{1}) \dot{\mu}_{00}}{\left(\ddot{x}_{1}\right)^3 x_{1}}}\right) \dot{\omega} \wedge \omega
-  d\mu_{00} \wedge \omega
\right]\wedge
$$
which admits the solution:
$$
\mu_{00}=F_2(\dot{x}_{1}x_{1}-2 x_{2})
$$
where $F_2(\dot{x}_{1}x_{1}-2 x_{2})$ is an arbitrary meromorphic function of $\dot{x}_{1}x_{1}-2 x_{2}$.

According to step 4, we must solve $\dg{M}=-M\mu$, where $M$ is scalar. $M$ is an at least first degree polynomial since, according to the lower bound given in the algorithm, we must have $\ord{M}\geq \frac{\ord{\mu}}{m}=1$.

Note that the equation $\dg{M}=-M\mu$ always admits the zero solution, which is not admissible. In addition it is readily seen that every other solution of degree greater than or equal to 1 with respect to $\ddt$ cannot be unimodular. Therefore the system is not flat.

Let us insist on the fact that a solution to the generalized moving frame structure equations is shown to exist for this non flat example, the only obstruction to flatness being the non existence of a unimodular integrating factor $M$.

\section{Concluding remarks}
We have proven that flatness is equivalent to the strong closedness of the ideal of 1-forms representing the differentials of all possible trivializations.
Moreover, we have separated the algebraic characterization of the differentials (Lemmas~\ref{Smith-lem} and \ref{Smith2-lem}) and the integrability aspects. The integrability conditions consist in finding solutions of a set of exterior differential equations which may be seen as a generalization, in the framework of manifolds of jets of infinite order, of the well-known moving frame structure equations. The computation of flat outputs in non trivial examples show the applicability of our results. Note that these flat outputs were already known since long, though obtained by \emph{ad hoc} methods. Other classes of examples, whose solutions are not presently known, are in preparation, using the computer algebra tools developed in \cite{AnL}.

We have chosen to present our results in terms of differential forms, in contrast with a large part of the control literature in this domain, where vector fields and Lie brackets are preferred. The main argument in favor of the former language is that the results are easy to state and  are constructive, without need to precise the coordinate system: flat outputs are directly computed by integration of the strongly closed ideal $\Omega$. 

Though the generalized moving frame structure equations are proved to be algebraically closed, a complete investigation of the existence of its local solutions in terms of involution remains to be done. In addition, there is no simple correlation between non flatness and non integrability of these equations: we have exhibited a non flat example for which such solutions exist, the only obstruction to flatness being that the unimodularity of the integrating factor cannot be satisfied.

Finally, note that many questions concerning upper bounds on polynomial degrees and on the number of derivatives of the state variables on which the flat outputs would depend still remain open.  In particular, our sequential procedure stops in a finite number of steps for flat systems, though this number is not known in advance, but no guarantee exists that the non flatness answer can also be obtained in a finite number of steps.

\paragraph{\textbf{Acknowledgements}} The author wishes to express his warm thanks to Michel Fliess, Philippe Martin, Pierre Rouchon and Philippe M\"{u}llhaupt for many helpful discussions and to Felix Antritter for his help in the computer algebraic aspects. He is particularly indebted to François Ollivier for his most fruitful discussions and suggestions. He is also grateful to Markku Nihtil\"a for his kind invitation to teach this topic at the University of Kuopio (Finland)\footnote{funded by a Marie Curie Host Fellowship for the Transfer of Knowledge, project PARAMCOSYS, MTKD-CT-2004-509223.}, and to Jose De Dona for his kind invitation to teach a similar course at the University of Newcastle (Australia).
\newpage

\appendix
\paragraph{\Large Appendix}
\section {The Smith decomposition algorithm}\label{Smith-app}
\label{polymat-sec}
We consider matrices of size $p\times q$, for  arbitrary integers $p$ and $q$, over the principal ideal domain $\kk[\ddt]$, here the non commutative ring of polynomials of $\ddt$ with coefficients in the field $\kk$ of meromorphic functions on a suitable time interval $\II$. The set of all such matrices is denoted by $\MM_{p,q}[\ddt]$.
For arbitrary $p\in \NN$, the set $\UU_{p}[\ddt]$ of unimodular matrices of size $p\times p$ is the subgroup of $\MM_{p,p}[\ddt]$ of invertible elements, namely the set of invertible polynomial matrices whose inverse is also polynomial. 

The following fundamental result on the
transformation of a polynomial matrix  over a principal ideal domain to its Smith form may be found in
\cite[Chap.8]{Co}):
\begin{thm}[Smith decomposition (or diagonal reduction)]\label{polymat-th-annex}
Given a $(\mu \times \nu)$ polynomial matrix $A$ over the non commutative ring $\kk[\ddt]$, there exist matrices
$V\in \UU_{\mu}[\ddt]$ and
$U\in \UU_{\nu}[\ddt]$
such that~$VAU=\left( \Delta, 0\right)$ if $\mu\leq \nu$ and~$VAU=\left( \begin{array}{c}\Delta\\ 0\end{array}\right)$ if $\mu\geq \nu$,
where
$\Delta$ is a $\mu\times \mu$ (resp. $\nu\times \nu$) diagonal matrix whose
diagonal elements, $(\delta_{1},\ldots,\delta_{\sigma},0,\ldots,0)$, are such
that
$\delta_{i}$ is a non zero $\ddt$-polynomial for $i=1,\ldots,\sigma$, and is a  divisor
of~$\delta_{j}$ for all~$\sigma\geq j\geq i$.
\end{thm}
The group of unimodular matrices admits a finite set of generators corresponding to the
following \emph{elementary right and left actions}:  
\begin{itemize}
\item \emph{right actions} consist of
permuting two columns, right multiplying a column by a non zero function of $\kk$, or adding the
$j$th column right multiplied by an arbitrary polynomial to the $i$th column, for
arbitrary $i$ and $j$; 
\item \emph{left actions} consist, analogously, of permuting two rows,
left multiplying a row by a non zero function of $\kk$, or adding the
$j$th row left multiplied by an arbitrary polynomial to the $i$th row, for
arbitrary $i$ and $j$.
\end{itemize}
Every elementary action may be represented by an \emph{elementary unimodular matrix} of the form $T_{i,j}(p)=I_{\nu}+1_{i,j}p$ with $1_{i,j}$ the matrix made of a single $1$ at the intersection of row $i$ and column $j$, $1\leq i,j\leq \nu$, and zeros elsewhere, with $p$ an arbitrary polynomial, and with $\nu=m$ for right actions and $\nu= n$ for left actions. 
One can easily prove that:
\begin{itemize}
\item right multiplication $AT_{i,j}(p)$ consists of adding the $i$th column  of $A$ right multiplied by $p$ to the $j$th column of A, the remaining part of $A$ remaining unchanged, 
\item left multiplication $T_{i,j}(p)A$ consists of adding the $j$th row of $A$ left multiplied by $p$ to the $i$th row of $A$, the remaining part of $A$ remaining unchanged,
\item $T_{i,j}^{-1}(p)=T_{i,j}(-p)$,
\item $T_{i,j}(1)T_{j,i}(-1)T_{i,j}(1)A$ (resp. $AT_{i,j}(1)T_{j,i}(-1)T_{i,j}(1)$) is the permutation matrix replacing the $j$th row of $A$ by the $i$th one and replacing the $j$th one of $A$ by the $i$th one multiplied by $-1$, all other rows remaining unchanged (resp. the permutation matrix replacing the $i$th column of $A$ by the $j$th one multiplied by $-1$ and replacing the $j$th one by the $i$th one, all other columns remaining unchanged).
\end{itemize}
Every unimodular matrix $V$ (left) and $U$ (right) may be obtained as
a product of such elementary unimodular matrices, possibly with a diagonal matrix $D(\alpha)=\diag{\alpha_{1},\ldots, \alpha_{\nu}}$ with $\alpha_{i}\in \kk$, $\alpha_{i}\neq 0$, $i=1,\ldots,\nu$, at the end since $T_{i,j}(p)D(\alpha)=D(\alpha)T_{i,j}(\frac{1}{\alpha_{i}}p\alpha_{j})$.

In addition, every unimodular matrix $U$ is obtained by such a product: its decomposition yields
$VU=I$ with $V$ finite product of the $T_{i,j}(p)$'s and a diagonal matrix. Thus, since the inverse of any $T_{i,j}(p)$ is of the same form, namely $T_{i,j}(-p)$, and since the inverse of a diagonal matrix is diagonal, it results that $V^{-1}=U$ is a product of elementary matrices of the same form, which proves the assertion.

The algorithm of decomposition of the matrix $A$ consists:
\begin{itemize}
\item first in permuting columns (resp. rows) to put the element of lowest degree in upper left
position, denoted by $a_{1,1}$, 
or creating this element by Euclidean division of two or more elements of the first row (resp. column) by suitable right actions (resp. left actions); 
\item then right divide all the other elements
$a_{1,k}$  (resp. left divide the $a_{k,1}$) of the new first row (resp. first column) by
$a_{1,1}$. If one of the rests is non zero, say
$r_{1,k}\neq 0$ (resp. $r_{k,1}\neq 0$), subtract the corresponding column (resp. row)
to the first column (resp. row) right multiplied (resp. left) by the corresponding quotient 
$q_{1,k}$ defined by the right
Euclidean division $a_{1,k}=a_{1,1}q_{1,k}+r_{1,k}$ (resp.
$q_{k,1}$ defined by $a_{k,1}=q_{k,1}a_{1,1}+r_{k,1}$). 
Then right multiplying all the columns by the corresponding quotients
$q_{1,k}$, $k=2,\ldots, \nu$ (resp. left multiplying rows by $q_{k,1}$, $k=2,\ldots,
\mu$), we iterate this process with the transformed first row (resp. first column) until it becomes $\left(a_{1,1},0,\ldots,0\right)$ (resp.  $\left(a_{1,1},0,\ldots,0\right)^{T}$ where $^{T}$ means transposition).
\item We then apply the same algorithm to the second row starting from $a_{2,2}$ and so on. To each transformation
of rows and columns correspond a left or right elementary unimodular matrix and
the unimodular matrix $V$ (resp. $U$) is finally obtained as the product of
all left (resp. right) elementary unimodular matrices so constructed.
\end{itemize}




\end{document}